\documentclass[12pt,reqno]{amsart}

\usepackage{amssymb, amsmath, amsfonts, latexsym}
\usepackage{enumerate}

\newtheorem{thm}{Theorem}[]
\newtheorem{lem}{Lemma}[section]
\newtheorem{cor}{Corollary}[section]
\newtheorem{prop}{Proposition}[]

\newtheorem{rmk}{Remark}[section]
\theoremstyle{definition}
\newtheorem{defn}{Definition}[section]

\numberwithin{equation}{section} \theoremstyle{remark}

\title[Kullback-Leibler distance]{\bf Approximation of Beta-Jacobi ensembles by Beta-Laguerre ensembles}

\author{Yutao Ma}
\address{Yutao MA\\ School of Mathematical Sciences $\&$ Laboratory  of Mathematics and Complex Systems of Ministry of Education, Beijing Normal University, 100875 Beijing, China.}
\thanks{The research of Yutao Ma was supported in part by NSFC 11571043, 11431014 and 985 Projects.}
\email{mayt@bnu.edu.cn}

\author{Xinmei SHEN}
\address{Xinmei SHEN\\ School of Mathematical Sciences, Dalian University of Technology, Linggong Road 2, Dalian, 116024, China.}
\thanks{The research of Xinmei Shen was supported in part by NSFC 11571058, HSSFMEC 17YJC910007 and FRF for the central universities DUT17LK31.}
\email{xshen@dlut.edu.cn}

\def\ml{\mathcal}
\def\f{\mathbf}

\def\<{\langle}
\def\>{\rangle}

\def\be{\begin{equation}}
\def\ee{\end{equation}}
\def\beaa{\begin{eqnarray*}}
\def\eeaa{\end{eqnarray*}}
\def\bea{\begin{eqnarray}}
\def\eea{\end{eqnarray}}
\def\lbl{\label}
\def\bdef{\begin{defn}}
\def\ndef{\end{defn}}

\def\bthm{\begin{thm}}
\def\nthm{\end{thm}}

\def\bprop{\begin{prop}}
\def\nprop{\end{prop}}

\def\brmk{\begin{remarks}}
\def\nrmk{\end{remarks}}

\def\bexa{\begin{exa}}
\def\nexa{\end{exa}}

\def\blem{\begin{lem}}
\def\nlem{\end{lem}}

\def\bcor{\begin{cor}}
\def\ncor{\end{cor}}

\date{}

\def\bexe{\begin{exe}}
\def\nexe{\end{exe}}

\def\bprf{\begin{proof}}
\def\nprf{\end{proof}}

\def\bdes{\begin{description}}
\def\ndes{\end{description}}

\begin{document}
\maketitle

\begin{abstract}
Let $\lambda$ and $\mu$ be beta-Jacobi and beta-Laguerre ensembles with joint density function $f_{\beta, m, a_1, a_2}$ and $f_{\beta, m, a_1}$, respectively. Here $\beta>0$ and $a_1, a_2$ and $m$ satisfying .  
$a_1, a_2>\frac{\beta}{2}(m-1).$ In this paper, we consider the distance between $2(a_1+a_2)\lambda$ and $\mu$ in terms of total variation distance and Kullback-Leibler distance.  Following the idea in \cite{JM2017}, we are able to prove that 
both the two distances go to zero once $a_1m=o(a_2)$ and not so if $\lim_{a_2\to\infty}a_1m/a_2=\sigma>0.$    
\end{abstract}



\section{Introduction}\label{chap:intro}
Let $\mu$ and $\nu$ be two probability measures on $(\mathbb{R}^n, \ml{B}),$
where $\mathbb{R}^n$ is the $n$-dimensional Euclidean space and $\ml{B}$ is the Borel $\sigma$-algebra.
We will consider the following two types distance between $\mu$ and $\nu:$

(1).  Total variation distance between $\mu$ and $\nu,$ denoted by $\|\mu-\nu\|_{\rm TV},$ is defined by
$$
\|\mu-\nu\|_{\rm TV}=2 \sup_{A\in \ml{B}}|\mu(A)-\nu(A)|=\int_{\mathbb{R}^n}|f(x)-g(x)|\, dx
$$
provided $\mu$ and $\nu$ have density functions $f$ and $g$ with respect to the Lebesgue measure, respectively.  

(2).  Kullback-Leibler distance between  $\mu$ and $\nu$ is defined by
\beaa\lbl{Kull}
D_{\rm KL}(\mu||\nu)=\int_{\mathbb{R}^n} \frac{d\mu}{d\nu}\log\frac{d\mu}{d\nu} d\nu.
\eeaa

Let $\beta>0$ be a constant  and $m\ge 1$ be an integer. A beta-Jacobi ensemble, also called the beta-MANOVA ensemble, is a set of random variables $\lambda:=(\lambda_1, \lambda_2, \cdots, \lambda_m)\in [0, 1]^m$ with joint probability density function
\be\lbl{bjdensity}
f_{\beta, a_1, a_2}(x_1, \cdots, x_m)=C_{\rm J}^{\beta, a_1, a_2}\prod_{1\le i<j\le m}|x_i-x_j|^{\beta}\prod_{i=1}^m x_i^{a_1-r}(1-x_i)^{a_2-r},
\ee
where $a_1, a_2>\frac{\beta}{2}(m-1)$ and $r:=1+\frac{\beta}{2}(m-1),$ and
$$C_{\rm J}^{\beta, a_1, a_2}=\prod_{j=1}^m \dfrac{\Gamma(1+\beta/2)\Gamma(a_1+a_2-\beta(m-j)/2)}{\Gamma(1+\beta j/2)\Gamma(a_1-\beta (m-j)/2)\Gamma(a_2-\beta(m-j)/2)}.$$

The density has close connections to the multivariate analysis of variance (MANOVA). For $\beta=1, 2, 4,$ the density function $f_{\beta, a_1, a_2}$ in \eqref{bjdensity} is the joint probability density function of the eigenvalues of independent matrices $\mathbf{Y}'\f Y(\f Y'\f Y+\f Z'\f Z)^{-1}$ with $a_1=\beta n_1/2$ and $a_2=\beta n_2/2.$ Here $\f Y=\f Y_{n_1\times m}$ and $\f Z=\f Z_{n_2\times m }$ are independent matrices with $n_1, n_2\ge m$ and the entries of both matrices are independent random variables with the standard real, complex or quaternion Gaussian distributions.  See \cite{Constantine} for $\beta=1$ and \cite{Muirhead} for $\beta=2,$ respectively.

A beta-Laguerre ensemble is a set of non-negative random variables $\lambda:=(\lambda_1, \lambda_2, \cdots, \lambda_m)$ with joint density function
\be\lbl{bldensity}
f_{\beta, \bar{a}}(x_1, \cdots, x_m)=C_{\rm L}^{\beta, \bar{a}}\prod_{1\le i<j\le m}|x_i-x_j|^{\beta}\prod_{i=1}^m x_i^{\bar{a}-r}e^{-\frac12\sum_{i=1}^m x_i},
\ee
where $\bar{a}>\frac{\beta}2(m-1)$ and $r=1+\frac{\beta}{2}(m-1),$ and
$$C_{\rm L}^{\beta, \bar{a}}=2^{-m\bar{a}}\prod_{j=1}^m \dfrac{\Gamma(1+\beta/2)}{\Gamma(1+(\beta/2)j)\Gamma(\bar{a}-(\beta/2)(m-j))}.$$
It is clear that $$\frac{f_{\beta, a_1, a_2}}{f_{\beta, a_1}}(x_1, \cdots, x_m)=\frac{C_{\rm J}^{\beta, a_1, a_2}}{C^{\beta, a_1}_{\rm L}} e^{\frac12\sum_{i=1}^m x_i} \prod_{i=1}^m (1-x_i)^{a_2-r}.$$

Let $\bold{\Gamma}_n=(\gamma_{ij})$ be a random orthogonal matrix which is uniformly distributed on
the orthogonal group $O(n).$ Let $\f Z_n$ be the $p_n\times q_n$ upper-left block of
$\bold{\Gamma}_n,$ where $p_n$ and $q_n$ are two positive integers.  
Denoted by $\ml{L}(\sqrt{n}\f Z_n)$  the joint probability distribution of the $p_nq_n$ random entries of $\sqrt{n}\f Z_n$
and ${\bf G}_n$ the joint distribution of $p_nq_n$ independent standard normals.  Let $f_n$ and $g_n$ be the probability density function of $\mathcal{L}(\sqrt{n}Z_n)$ and $\mathcal{L}({\bf G}_n)$ with respect to the Lebesgue measure, respectively.  
According to the explicit expression of  $f_n/g_n$ in \cite{JM2017}, it has a particular form of $\dfrac{f_{\beta, a_1, a_2}}{f_{\beta, a_1}}$ with $\beta=1, \; m=q, \; a_1=\frac{p}{2}$ and $a_2=\frac{n-p}{2}.$ 
In \cite{Jiang06}, 
Jiang proves that when $p_n=o(\sqrt{n})$ and $q_n=o(\sqrt{n})$,
$$\lim_{n\to\infty}\|\ml{L}(\sqrt{n}\f Z_n)-{\bf G}_n\|_{\rm TV}=0$$ while
when $p_n=O(\sqrt{n})$ and $q_n=O(\sqrt{n}),$
$$\liminf_{n\to\infty}\|\ml{L}(\sqrt{n}\f Z_n)-{\bf G}_n\|_{\rm TV}>0.$$
This is the first result to characterize exactly how many entries of a typical orthogonal matrix could be approximated by independent standard normals. 
Recently, Jiang and the first author in \cite{JM2017} completely resolve this problem. Precisely, they show that 
$$\aligned  \lim_{n\to\infty}d(\ml{L}(\sqrt{n}\f Z_n), {\ml L}({\bf G}_n))=0, &\quad {\rm if } \quad pq=o(n);  \\
\liminf_{n\to\infty}d(\ml{L}(\sqrt{n}\f Z_n), {\ml L}({\bf G}_n))>0,  &\quad {\rm if} \quad pq=O(n). \endaligned $$
Here $d$ is the total variation distance, Kullback-Leibler distance or  Hellinger distance.  
In 2013, Jiang in \cite{Jiang13} works on general $\beta>0.$   He proves that when
\be\lbl{condition}\aligned
 m &\to \infty, \quad a_1\to\infty \quad {\rm and } \quad a_2\to \infty \quad {\rm such \;  \; that } \; \\
a_1 &=o(\sqrt{a_2}), \quad m=o(\sqrt{a_2}) \quad\mbox{and}\;  \; \frac{m\beta}{2a_1}\rightarrow \gamma\in(0,1],
\endaligned \ee
it holds
$$\lim_{a_2\to\infty}\|\mathcal{L}(2a_2 \lambda)-\mathcal{L}(\mu)\|_{\rm TV}=0,$$
where $\lambda=(\lambda_1, \cdots, \lambda_m)$ have joint probability density function $f_{\beta, a_1, a_2}$ as in \eqref{bjdensity} and $\mu=(\mu_1, \cdots,\mu_m)$  have joint probability density function $f_{\beta, a_1}$ as in \eqref{bldensity}.  

Inspired by the work in \cite{Jiang13}  and \cite{JM2017}, for general $\beta>0,$ we want to completely understand the behavior between $\lambda$ and $\mu.$ 
Making a minor adjustment from $d\big({\ml L}(2a_2\lambda), {\ml L}(\mu)\big)$ in \cite{Jiang13}, we will investigate the following object
$$d\big({\ml L}(2a\lambda), {\ml L}(\mu)\big)$$ under the condition $a_1m=o(a_2)$ or $a_1m=O(a_2)$ with $a:=a_1+a_2.$     

For two different distances mentioned above, we have the following theorem. 
\bthm\label{bjconv} Let $\mu=(\mu_1, \mu_2, \cdots, \mu_m)$ and $\lambda=(\lambda_1, \lambda_2, \cdots, \lambda_m)$ be random variables with density $f_{\beta, a_1}$ as in \eqref{bldensity} and
$f_{\beta, a_1, a_2}$ as in \eqref{bjdensity}, respectively.  Let $d\big(\ml{L}(2a\lambda), \ml{L}(\mu)\big)$ be the total variation distance or the Kullback-Leibler  distance
 between the probability distributions of $2a\lambda$ and $\mu.$ 
Then 
\begin{itemize}
\item[(i).]  $\lim_{a_2\to\infty} d\big(\ml{L}(2a\lambda), \ml{L}(\mu)\big)=0$ if $a_1m=o(a_2).$
\item[(ii).]  $\liminf_{a_2\to\infty}d\big(\mathcal{L}(2a \lambda), \mathcal{L}(\mu)\big)>0$ if $\lim_{a_2\to\infty}\frac{a_1m}{a_2}=\sigma>0.$
\end{itemize}
\nthm
By Pinsker's inequality, we know
\be\lbl{TV_KL}
 \|\mu-\nu\|_{\rm TV}^2\le 2 D_{\rm KL}(\mu||\nu). 
\ee
Therefore as in \cite{JM2017}, for the first item, we just need to prove 
\be\lbl{kl=0}\lim_{a_2\to\infty}D_{\rm KL}\big(\mathcal{L}(2a \lambda)||\mathcal{L}(\mu)\big)=0\ee
and for the second item it suffices to prove 
\be\lbl{tv>0}\liminf_{a_2\to\infty}  \|\ml{L}(2a\lambda)-\ml{L}(\mu)\|_{\rm TV}>0.\ee
Furthermore, for the validity of \eqref{tv>0}, by Lemma 2.15 in \cite{JM2017}, it is enough to prove \eqref{tv>0} under  the following three conditions:
\begin{itemize}
\item[\bf A1.] $m\equiv 1$ and $ \lim_{a_2\to\infty} \frac{a_1}{a_2}\in (0, 1);$ 
\item[\bf A2.] $m\to\infty, \; \lim_{a_2\to\infty}\frac{m}{a_1}=0$ and $\lim_{a_2\to\infty}\frac{ma_1}{a_2}=\sigma>0;$
\item[\bf A3.]  $m\to\infty, \; \lim_{a_2\to\infty}\frac{a_1}{\sqrt{a_2}}=x $ and $ \lim_{a_2\to\infty}\frac{m}{\sqrt{a_2}}=y.$ 
\end{itemize} 

Set $\eta=\frac{\beta}2$ and define
$$\aligned
& K_m=(\frac{1}{a})^{m a_1}\prod_{i=0}^{m-1} \frac{\Gamma(a-\eta i)}{\Gamma(a_2-\eta i)}, \\
& L_m(x_1, \cdots, x_m)=e^{\frac12\sum_{i=1}^m x_i}\prod_{i=1}^{m}(1-\frac{x_i}{2a})^{a_2-r} {\f I}_{\{\max x_i\le 2a\}}.
\endaligned$$

We will show in the forth section that the total variation distance could be regarded as 
\be\lbl{TV-e}
 \|\ml{L}(2a\lambda)-\ml{L}(\mu)\|_{\rm TV}=\mathbb{E} |K_m L_m(\mu)-1|. \ee  
Meanwhile for the Kullback-Leibler distance, 
 we understand it as  
\be\lbl{KL-e} D_{\rm KL}\big(\mathcal{L}(2a \lambda)||\mathcal{L}(\mu)\big)=\mathbb{E} \log \big(K_m L_m(\lambda)\big).
\ee
To prove \eqref{kl=0}, with the help of the expression \eqref{KL-e} and Taylor's formula for $\log L_m$, one just needs to characterize the asymptotics of $\log K_m$ and  
to have the asymptotical expression for $\sum_{i=1}^m\mathbb{E}\lambda_i^{k}$ with $k=1,2,3,$ where $(\lambda_1, \cdots, \lambda_m)$ have joint density function $f_{\beta, a_1, a_2}.$ 
According to the interpretation of Edelman and Sutton in \cite{ES} (see also \cite{DP}),  $f_{\beta, a_1, a_2}$ is also the joint density function of  the eigenvalues of $\f B\f B'.$ The explicit form of
$m$ by $m$ random matrix $\f B$ will be given later in \eqref{matrixb}, whose elements are related to mutually independent Beta distributions. There isn't any result on  $\sum_{i=1}^m\mathbb{E}\lambda_i^{k}$ when $m a_1=o(a_2),$ 
which then requires tendious calculations related to Beta distribution presented in Section 2.    

The proof of \eqref{tv>0} is much more difficult. We have to establish a central limit theorem for $\log (K_m L_m(\mu))$ as in \cite{JM2017}.  Review $r=1+\frac{\beta(m-1)}{2}$ and set \be\lbl{defu} 
U_m:=\frac{r}{2a_2}\sum_{i=1}^m (\mu_i-2a_1)-\frac{(a_2-r)}{8a_2^2}\sum_{i=1}^m(\mu_i-2a_1)^2. 
\ee 
With the help of Taylor's formula and the property of logarithmic Gamma function, we are able to write 
$$\log(K_mL_m(\mu))=U_m-\mathbb{E}U_m+C_m.$$ Here $C_m$ converges to some constant in probability as $a_2\to\infty$ when either {\bf A2} or {\bf A3} is satisfied. 
Therefore, the key task for us is to obtain the central limit theorem for $U_m-\mathbb{E}U_m$ as follows. 
\bprop\lbl{uniqprop}  Let $\mu=(\mu_1, \mu_2, \cdots, \mu_m)$ be random variables with density $f_{\beta, a_1}$ as in \eqref{bldensity} and $U_m$ be given by \eqref{defu}. 
Then under the assumption {\bf A2} or {\bf A3}, with $\sigma:=xy$ in {\bf A3} we have 
$$U_m+\frac{(a_2-r)a_1mr}{2a_2^2}\to N(0, \frac{\beta \sigma^2}{4})$$
weakly as $a_2\to\infty.$ Here $\frac{(a_2-r)a_1mr}{2a_2^2}=-\mathbb{E}U_m.$  
\nprop  

\begin{rmk} Theorem 1.5 in \cite{DE06} tells that $U_m-c_m$ converges weakly to some normal distribution under the assumption {\bf A3} as $a_2\to\infty$. Here $c_m=d_m\mathbb{E}U_m$ with $d_m$ is given as 
$$d_m=\frac{2a_2}{\gamma(a_2-r)}(\gamma-\frac{\beta m}{2a_1})+\frac{a_1}{r}\big(1-\frac{\beta m}{\gamma a_1}+\frac{m^2\beta^2(1+\gamma)}{4a_1^2\gamma^2}\big).$$
It is easy to check that $$\lim_{a_2\to\infty}d_m=1\quad\text{and}\quad (d_m-1)r=o(a_1)$$ under the assumption  {\bf A3} as $a_1$ large enough.  
Obviously, $\lim_{a_2\to\infty} \frac{\mathbb{E}U_m}{c_m}=1.$
 However, since $\mathbb{E}U_m$ has the same order as $r$ under the assumption  {\bf A3}, we know   
$$c_m-\mathbb{E}U_m=(d_m-1)\mathbb{E}U_m=o(a_1).$$  Then $c_m-\mathbb{E}U_m$ is not necessarily finite under the assumption  {\bf A3}, which claims the failure of the
central limit Theorem for $U_m-\mathbb{E} U_m$ via Theorem 1.5 in \cite{DE06}.  This cruel fact forces us to work very hard directly on the central limit theorem for $U_m-\mathbb{E}U_m.$  
\end{rmk}  
For this aim, we need the characterization of Dumitriu and Edelman in their famous work \cite{DE02}.  They understand $f_{\beta, a_1}$ as the joint density function of the eigenvalues of  the random matrix $\f A\f A^{\prime}.$ The $m$ by $m$ random matrix $\f A$ will be introduced later in \eqref{matrixa}, whose elements are mutually independent chi distribution.  Based on this characterization, by independence and the properties of chi square distribution, in the third section we prove Proposition \ref{uniqprop} via the central limit theorem for the sum of independent random variables under {\bf A2} and that for $m$-dependent random variables under {\bf A3}, respectively.     

Therefore this paper will be organized as follows:

\noindent\textbf{Section \ref{prepara}: Preliminaries}  

Section \ref{asymK}:  On the asymptotics of $K_m.$ 

Section \ref{triplepro}: On $(\lambda_1, \lambda_2, \cdots, \lambda_m)$ having joint density function $f_{\beta, a_1, a_2}.$

Section \ref{bipro}:  On $(\mu_1, \mu_2, \cdots, \mu_m)$ having joint density function $f_{\beta, a_1}.$

\noindent\textbf{Section \ref{prfprop}: Proof of  Proposition \ref{uniqprop}}

Section \ref{centralA2}: The proof of Proposition \ref{uniqprop} under {\bf A2}.

Section \ref{centralA3}: The proof of Proposition \ref{uniqprop} under {\bf A3}. 

\noindent \textbf{Section \ref{prfthm}:  Proof of Theorem \ref{bjconv}}. 

\section{Preliminaries}\lbl{prepara}
In this section, we collect all Lemmas and Propositions we need. 
\subsection{On the asymptotics  of $K_m$.}\lbl{asymK}
\blem\lbl{kmO} For $0\le \beta (m-1)< 2a_1<2a,$ recall 
\be\lbl{K}
K_m=(\frac{1}{a})^{m a_1}\prod_{i=0}^{m-1} \frac{\Gamma(a-\eta i)}{\Gamma(a_2-i\eta)}.\ee  Suppose $a_1\to\infty, \, \limsup_{a_2\to\infty}\frac{a_1}{a}<1$  and $a_1m=O(a_2)$ as $a_2\to\infty.$  Then 
$$\aligned \log K_m=-a_1m+m\big(a_2-\frac{r}{2}\big)\log(1+\frac{a_1}{a_2})-\frac{\beta^2 a_1m^3}{24a^2}+o(1).
\endaligned $$
\nlem

\bprf Recall Stirling's formula:
$$\log\Gamma(x)=(x-\frac12)\log x-x+\log\sqrt{2\pi}+\frac1{12x}+O(\frac{1}{x^3})$$
as $x\to+\infty.$
Therefore applying Stirling's formula to $\log\Gamma(a-\eta i)$ and $\log\Gamma(a_2-\eta i)$ and combining alike terms,  we have 
$$\aligned \log K_m&=-ma_1\log a+\sum_{i=0}^{m-1}\log\Gamma(a-\eta i)-\sum_{i=0}^{m-1}\log\Gamma(a_2-\eta i)\\
&=-ma_1\log a+\sum_{i=0}^{m-1}\Big[(a-\eta i-\frac12)\log(a-\eta i)\\
&\quad\quad\quad\quad\quad\quad\quad -(a_2-\eta i-\frac12)\log(a_2-\eta i)-a_1\Big]
+o(1).\\
\endaligned
$$
By writing $a-\eta i=a_2-\eta i+a_1$ and putting the term $-ma_1\log a$  into the sum $\sum_{j=0}^{m-1},$ we see 
$$\log K_m
=-a_1m+\sum_{i=0}^{m-1}(a_2-\eta i-\frac12)\log(\frac{a-\eta i}{a_2-\eta i})+a_1\sum_{i=0}^{m-1}\log(1-\frac{\eta i}{a})+o(1).$$
Applying the decomposition
$$\log(\frac{a-\eta i}{a_2-\eta i})=\log(1+\frac{a_1}{a_2})+\log(1+\frac{\eta a_1 i}{a(a_2-\eta i)})$$
to the expression of $\log K_m$ above and by the fact $$\sum_{i=0}^{m-1} (a_2-\eta i-\frac12)=(a_2-\frac{1+\eta(m-1)}{2})m=(a_2-\frac{r}2)m,$$  we have  
\be\lbl{Kexpres}\aligned
\log K_m&=-a_1m+m\big(a_2-\frac{r}{2}\big)\log(1+\frac{a_1}{a_2})+o(1)\\
&\quad +\sum_{i=0}^{m-1}\big((a_2-\eta i-\frac12)\log\big(1+\frac{\eta a_1 i}{a(a_2-\eta i)}\big)+a_1\log(1-\frac{\eta i}{a})\big).\endaligned \ee
Since $\log(1+x)=x-\frac{x^2} 2+O(x^3)$ and $\log(1+x)=x+O(x^2)$ as $x\to 0,$  we have 
\be\lbl{expres1}\aligned 
&a_1\log(1-\frac{\eta i}{a})=-\frac{\eta a_1 i}{a}-\frac{a_1\eta^2 i^2}{2a^2}+O(\frac{a_1m^3}{a^3})\\
&(a_2-\eta i)\log\big(1+\frac{\eta a_1 i}{a(a_2-\eta i)}\big)=\frac{\eta a_1 i}{a}+
O(\frac{a_1^2m^2}{a_2^3})\\
&-\frac12\log\big(1+\frac{\eta a_1 i}{a(a_2-\eta i)}\big)=O(\frac{a_1m}{a_2^2})
\endaligned \ee
for any $0\le i\le m-1.$
The condition $\beta (m-1)< 2a_1$ and $a_1m=O(a_2)$ implies that all the following three terms
$O(\frac{a_1m^3}{a^3}), \;O(\frac{a_1^2m^2}{a_2^3})$ and $O(\frac{a_1m}{a_2^2})$
 could be written as $o(\frac1m).$ Therefore  it follows from \eqref{expres1} that  
\be\lbl{sumexpr}\aligned &\quad \sum_{i=0}^{m-1}\big(a_1\log(1-\frac{\eta i}{a})+(a_2-\eta i-\frac12)\log\big(1+\frac{\eta a_1 i}{a(a_2-\eta i)}\big)\big) \\
&=\sum_{i=0}^{m-1}-\frac{a_1\eta^2 i^2}{2a^2}+o(1)\\
&=-\frac{\eta^2 a_1m^3}{6a^2}+o(1).\endaligned  \ee
Now putting \eqref{sumexpr} back into \eqref{Kexpres}, we have  $$\log K_m=-a_1m+m\big(a_2-\frac{r}{2}\big)\log(1+\frac{a_1}{a_2})-\frac{\beta^2 a_1m^3}{24a^2}+o(1).
$$
The proof is complete. 
\nprf

\subsection{On $(\lambda_1, \lambda_2, \cdots, \lambda_m)$ having joint density function $f_{\beta, a_1, a_2}$ }\lbl{triplepro} 
Now we want to understand what $\mathbb{E}\sum_{i=1}^m \lambda_i^k$  will be alike for $k=1, 2, 3$ when $a_1 m=o(a_2).$ 
However this asymptotic could not be provided by  the explicit form of the joint density \eqref{bjdensity}.  Therefore we need the help of
the interpretation from Edelman and Sutton \cite{ES} (see also \cite{DP}) as mentioned in the Introduction.  That is,  the eigenvalues of $\f B\f B'$ have joint density function $f_{\beta, a_1, a_2},$
where the $m$ by $m$ random matrix $\f B$ has the form  \begin{align}\label{matrixb}
\f B=
    \begin{pmatrix}
      \sqrt{c_ms'_{m-1}}&  &  &   &  \\
      -\sqrt{s_{m-1}c'_{m-1}}& \sqrt{c_{m-1}s'_{m-2}} &   &   &\\
       &  -\sqrt{s_{m-2}c'_{m-2}} & \sqrt{c_{m-2} s'_{m-3}} &   &\\
      &    &   \ddots &\ddots & \\
      &   &   &   -\sqrt{s_1c'_1} & \sqrt{c_1}
    \end{pmatrix}
\end{align}
with the non-negative random variables $c_i, s_i$ with $1\le i\le m$ and $c'_i,  s'_i$ with $1\le  i\le m-1$ obeying the distribution and relationships
\begin{itemize}
\item[1).] $\{c_1, c_2, \cdots, c_m, c'_1, c'_2, \cdots, c'_{m-1}\}$ mutually independent;
\item[2).] $c_i\sim {\rm Beta}(a_1-\eta(m-i), a_2-\eta(m-i));$
\item[3).] $c'_i\sim{\rm Beta}(\eta i, a_1+a_2-\eta(2m-i-1));$
\item[4).] $s_i+c_i=1, \quad s'_i+c'_i=1.$
\end{itemize}

Based on this interpretation, Dumitriu and Paquette \cite{DP} obtained a series expansion of the scaled moment $\frac{1}{m}\mathbb{E}{\rm tr}((\f B\f B')^k)$ when $a_1, m$ and $a_2$ have same order.  Precisely,
$
\frac{1}{m}\mathbb{E}{\rm tr}((\f B\f B')^k)=\sum_{j=0}^\infty \rho_k(j,\alpha)m^{-j}.
$ 
The coefficients $\rho_k(j,\alpha)$ are palindromic polynomials in $(-\alpha)$ of degree $j$. This result is perfect with concise form. However, it is too hard to have a direct form via this characterization and 
it does not satisfy the assumption neither {\bf A2} nor  {\bf A3}. Therefore, with the help of \eqref{matrixb}, we calculate directly the following expressions under $a_1m=o(a_2).$

\bprop\lbl{objsquare} Suppose that $a_1m=o(a_2)$ as $a_2\to\infty.$  Assume that $(\lambda_1, \lambda_2, \cdots, \lambda_m)$ have joint probability density function $f_{\beta, a_1, a_2}$ as in \eqref{bjdensity}.  Then we have 
\be\lbl{e0}
\aligned
\mathbb{E}\sum_{i=1}^{m}\lambda_i&=\frac{a_1m}{a}+o(ma_2^{-1});\\
\mathbb{E}\sum_{i=1}^m\lambda_i^2&=\frac{1}{a^2}(a_1^2m+\eta  a_1 m^2)+o(a_2^{-1});\\
\mathbb{E}\sum_{i=1}^m\lambda_i^3&=\frac{1}{a^3}\big( a_1^3m+3\eta a_1^2m^2+\eta^2a_1m^3\big)+o(a_2^{-1})\\
\endaligned
\ee
as $a_2\to \infty$ with $a=a_1+a_2,$

\nprop
\bprf By the interpretation above, with the convention $s_0'=1,$ we have
\be\label{eo}
\mathbb{E}\sum_{i=1}^m\lambda_i=\mathbb{E}{\rm tr}(\f B\f B')=\sum_{i=1}^m\mathbb{E}c_{m+1-i} s'_{m-i}+\sum_{i=1}^{m-1}\mathbb{E}s_{m-i}c'_{m-i}\ee 
and 
\be\label{e}
\aligned
\mathbb{E}\sum_{i=1}^{m}\lambda_i^2&=\mathbb{E} {\rm tr}((\f B\f B')^2)\\
&=\mathbb{E}\sum_{i=1}^{m-1}s_{m-i}^2 (c'_{m-i})^2+\mathbb{E}\sum_{i=1}^{m}c_{m+1-i}^2(s'_{m-i})^2\\
&\quad+2\mathbb{E}\sum_{i=1}^{m-1} s_{m-i}c'_{m-i}c_{m-i}s'_{m-i-1}+2\mathbb{E}\sum_{i=1}^{m-1}c_{m+1-i}c'_{m-i}s_{m-i}s'_{m-i}.
\endaligned \ee
According to the expressions \eqref{eo} and \eqref{e}, for $\mathbb{E}\sum_{i=1}^m\lambda_i$ and $\mathbb{E}\sum_{i=1}^{m}\lambda_i^2,$ we have to work on  the following six items:
$$\mathbb{E}c_{m-i}, \quad \mathbb{E}c'_{m-i}, \quad \mathbb{E}c_{m-i}^2, \quad \mathbb{E}(c'_{m-i})^2, \quad \mathbb{E}s_{m-i}^2 \quad\text{and}\quad \mathbb{E}(s'_{m-i})^2.$$
For the random variable $\xi\sim {\rm Beta}(\alpha, \beta),$ it is well-known that 
\be\lbl{betavar}\mathbb{E} \xi=\frac{\alpha}{\alpha+\beta} \quad {\rm and}\quad \mathbb{E}\xi^2=\frac{\alpha (\alpha+1)}{(\alpha+\beta)(\alpha+\beta+1)}.
\ee

Since $2a_1> \beta(m-1),$ it enforces that $m^2/a_2\to 0$ when $a_2\to\infty$ from $a_1m=o(a_2).$   
By definition and \eqref{betavar}, keeping in mind that $a=a_1+a_2,$ one gets
\be\label{fouro}
\mathbb{E}c_{m-i}=\frac{a_1-\eta i}{a-2\eta i}=\frac{a_1-\eta i}{a}+o(a_2^{-1})
\ee 
for $0\le i\le m-1.$
Here and after we use frequently the following trick to make $i$ vanish from the denominator as for \eqref{fouro}. That is
 $$\frac{a_1-\eta i}{a-2\eta i}=\frac{a_1-\eta i}{a}(\frac{a-2\eta i+2\eta i}{a-2\eta i})=\frac{a_1-\eta i}{a}+\frac{2\eta i(a_1-\eta i)}{a(a-2\eta i)}=\frac{a_1-\eta i}{a}+o(a_2^{-1}),$$
 where the last equality holds since $i(a_1-\eta i)\le a_1 m=o(a_2)=o(a).$ Similarly, we have 
\be\lbl{four1}\aligned 
&\mathbb{E}c'_{m-i}=\frac{\eta(m-i)}{a-\eta(2i-1)}=\frac{\eta(m-i)}{a}+o(a_2^{-1}); \\
&\mathbb{E}c_{m-i}^2=\frac{(a_1-\eta i)(a_1-\eta i+1)}{(a-2\eta i)(a-2\eta i+1)}=\frac{(a_1-\eta i)^2+(a_1-\eta i)}{a^2}+o(\frac{a_1}{a^2});\\
&\mathbb{E}(c'_{m-i})^2=\frac{\eta(m-i)(\eta(m-i)+1)}{(a-2\eta i+\eta)(a-2\eta i+\eta+1)}=\frac{\eta^2(m-i)^2}{a^2}+o(a_2^{-3/2}).
 \endaligned
 \ee
Consequently
\be\lbl{four2}\aligned &\mathbb{E} s_{m-i}^2=\mathbb{E} (1-c_{m-i})^2=1-\frac{2(a_1-\eta i)}{a}+\frac{(a_1-\eta i)^2+(a_1-\eta i)}{a^2}+o(a_2^{-1});\\
&\mathbb{E}( s'_{m-i})^2=\mathbb{E} (1-c'_{m-i})^2=1-\frac{2\eta(m-i)}{a}+o(a_2^{-1}).
\endaligned \ee 

Plugging \eqref{fouro} and \eqref{four1} into the expression \eqref{eo},  we have
 $$\aligned \mathbb{E}\sum_{i=1}^m\lambda_i&=\sum_{i=1}^m\mathbb{E}c_{m+1-i}+\sum_{i=1}^{m-1}\mathbb{E}c_{m-i}^{\prime}-\sum_{i=1}^m\mathbb{E}c_{m+1-i}\mathbb{E}c_{m-i}^{\prime}-\sum_{i=1}^{m-1}\mathbb{E}c_{m-i}\mathbb{E}c_{m-i}^{\prime}\\
 &=\frac{1}{a}\sum_{i=1}^{m-1}\bigg(a_1-\eta(i-1)+\eta (m-i)\bigg)+\frac{a_1}{a}+o(ma_2^{-1})\\
 &=\frac{a_1m}{a}+o(ma_2^{-1}).\endaligned
  $$
Here for the second term, we use the facts $a_1m=o(a_2),$ \eqref{fouro} and \eqref{four1} to get $$\sum_{i=1}^m\mathbb{E}c_{m+1-i}\mathbb{E}c_{m-i}^{\prime}=\sum_{i=1}^{m}O(a_1ma^{-2})=O(a_1m^2 a_2^{-2})=o(ma_2^{-1}).$$
Next we focus on the second expression in \eqref{e0}. We treat the first term of  \eqref{e}.
Since $m=o(\sqrt{a_2})$ and  $o(ma_2^{-3/2})=o(a_2^{-1}),$ we can drop off the terms of order $o(a_2^{-s})$ with $s\ge 3/2$ in the sum $\sum_{i=1}^m.$ This would greatly simplify 
the calculus.  
Thereby, based on \eqref{four1}, \eqref{four2} and the condition $a_1m=o(a_2)$ and $\beta (m-1)< 2a_1,$ it follows from the independence of 
$\{c_1, \cdots, c_m, c_1', \cdots, c'_{m-1}\}$ that 
$$\aligned \mathbb{E}\sum_{i=1}^{m-1}s_{m-i}^2 (c'_{m-i})^2&=\sum_{i=1}^{m-1}\big((1-\frac{a_1-\eta i}{a})^2+\frac{a_1-\eta i}{a^2}+o(a_2^{-1})\big)\big(\frac{\eta^2(m-i)^2}{a^2}+o(a_2^{-3/2})\big)\\
&=\sum_{i=1}^{m-1}\big(\frac{\eta^2(m-i)^2}{a^2}+o(a_2^{-3/2})\big)\\
&=\frac{\eta^2 m^3}{3a^2}+o(a_2^{-1}),
\endaligned $$
where for the second equality we drop off directly the term 
$$\bigg(-\frac{2(a_1-\eta i)}{a}+\frac{(a_1-\eta i)^2}{a^2}+\frac{a_1-\eta i}{a^2}+o(a_2^{-1})\bigg)\frac{\eta^2(m-i)^2}{a^2}=o(a_2^{-3/2})$$
 for any $1\le i\le m-1.$ 
Similarly since $a_1 m=o(a_2),$ we have 
$$\aligned \mathbb{E}\sum_{i=1}^{m}c_{m+1-i}^2(s'_{m-i})^2&=\sum_{i=1}^m \big(\frac{(a_1-\eta i+\eta)^2+(a_1-\eta i+\eta)}{a^2}+o(\frac{a_1}{a^2})\big)\big(1-\frac{2\eta(m-i)}{a}+o(a_2^{-1})\big)\\
&=\sum_{i=1}^m\big(\frac{(a_1-\eta i)^2}{a^2}+O(\frac{a_1}{a^2})\big)\\
&=\frac{3a_1^2m-3\eta a_1m^2+\eta^2m^3}{3a^2}+o(a_2^{-1}).
\endaligned $$
The same argument also leads  
$$\aligned \mathbb{E}\sum_{i=1}^{m-1} s_{m-i}c'_{m-i}c_{m-i}s'_{m-i-1}&=\sum_{i=1}^{m-1}(\mathbb{E}c_{m-i}-\mathbb{E}c^2_{m-i})\mathbb{E}c'_{m-i}(1-\mathbb{E}c'_{m-i-1})
\\
&=\sum_{i=1}^{m-1}\bigg(\frac{\eta(a_1-\eta i)(m-i)}{a^2}+o(\frac{a_1}{a^2})\bigg)\\
&=\frac{3\eta a_1m^2-\eta^2 m^3}{6a^2}+o(a_2^{-1})
\endaligned $$
and 
$$\aligned 
\mathbb{E}\sum_{i=1}^{m-1}c_{m+1-i}c'_{m-i}s_{m-i}s'_{m-i}&=\sum_{i=1}^{m-1}\bigg(\frac{\eta(a_1-\eta i)(m-i)}{a^2}+o(\frac{a_1}{a^2})\bigg)\\
&=\frac{3\eta a_1m^2-\eta^2 m^3}{6a^2}+o(a_2^{-1}). \endaligned $$
Therefore plugging all these four expressions above into \eqref{e}, we have 
 $$\aligned 
\mathbb{E}\sum_{i=1}^{m}\lambda_i^2&=\frac{\eta^2 m^3}{3a^2}+\frac{3a_1^2m-3\eta a_1m^2+\eta^2m^3}{3a^2}+4\cdot\frac{3\eta a_1m^2-\eta^2 m^3}{6a^2}+o(a_2^{-1})\\
&=\frac{a_1^2m+\eta a_1 m^2}{a^2}+o(a_2^{-1}).
\endaligned  $$ 
Now we work on the last expression in \eqref{e0}. 
For the Beta distribution  $\xi\sim {\rm Beta}(\alpha, \beta),$ one knows 
$$\mathbb{E}\xi^3=\frac{\alpha(\alpha+1)(\alpha+2)}{(\alpha+\beta)(\alpha+\beta+1)(\alpha+\beta+2)}.$$
 Therefore we have
 \be\label{xi3}\aligned
 &\mathbb{E}c_{m-i}^3 =\frac{(a_1-\eta i)^3}{a^3}+o(a_1a_2^{-2});\\
& \mathbb{E}(c'_{m-i})^3 =\frac{\eta^3(m-i)^3}{a^3}+o(a_1a_2^{-2}).
 \endaligned \ee
With careful calculation, we have
$$\aligned \mathbb{E}\sum_{i=1}^m \lambda_i^3&=\sum_{i=1}^m \mathbb{E}c_{m-i+1}^3(1-\mathbb{E}c'_{m-i})^3+\sum_{i=1}^{m-1}\mathbb{E}(c'_{m-i})^3(1-\mathbb{E}c_{m-i})^3\\
&\quad+3\sum_{i=1}^{m-1}(1-\mathbb{E}c_{m-i})^2\bigg(\mathbb{E}(c'_{m-i})^2-\mathbb{E}(c'_{m-i})^3\bigg)\mathbb{E}c_{m-i+1}\\
&\quad+3\sum_{i=1}^{m-1}(1-\mathbb{E}c'_{m-i-1})^2\bigg(\mathbb{E}c_{m-i}^2-\mathbb{E}c_{m-i}^3\bigg)\mathbb{E}c'_{m-i}\\
&\quad+3\sum_{i=1}^{m-1}(1-\mathbb{E}c_{m-i})\mathbb{E}\bigg(c'_{m-i}(1-c'_{m-i})^2\bigg)\mathbb{E}c_{m-i+1}^2\\
&\quad+3\sum_{i=1}^{m-1} (1-\mathbb{E}c'_{m-i-1})\mathbb{E}\bigg(c_{m-i}(1-c_{m-i})^2\bigg)\mathbb{E}(c'_{m-i})^2\\
&\quad+3\sum_{i=1}^{m-1}\mathbb{E}(c_{m-i}-c^2_{m-i}) \mathbb{E}\bigg(c'_{m-i}-(c'_{m-i})^2\bigg)\mathbb{E}c_{m-i+1}(1-\mathbb{E}c'_{m-i-1})\\
&\quad+3\sum_{i=1}^{m-1}\mathbb{E}(c_{m-i+1}-c^2_{m-i+1}) \mathbb{E}\bigg(c'_{m-i}-(c'_{m-i})^2\bigg)\mathbb{E}c'_{m-i+1} (1-\mathbb{E}c_{m-i}).
\endaligned $$
Similarly as for $\mathbb{E}\sum_{i=1}^m \lambda_i^2,$ we drop off the terms $o(a_1a_2^{-2}).$ Plugging  \eqref{fouro}, \eqref{four1} and \eqref{xi3} into above expression, one gets
 $$\aligned \mathbb{E}\sum_{i=1}^m \lambda_i^3&=\sum_{i=1}^m\frac{(a_1-\eta i)^3}{a^3}+\sum_{i=1}^{m-1}\frac{\eta^3(m-i)^3}{a^3}+o(a_2^{-1})\\
 &\quad+\frac{9}{a^3}\sum_{i=1}^{m-1}\big\{\eta(m-i)(a_1-\eta i)^2+\eta^2(m-i)^2(a_1-\eta i)\big\}\\
 &=\frac{1}{a^3}\big(a_1^3m+3\eta a_1^2m^2+\eta^2a_1m^3\big)+o(a_2^{-1}).\\
  \endaligned $$
  This finally closes the entire proof. 
\nprf 

\subsection{On $\mu=(\mu_1, \mu_2, \cdots, \mu_m)$ having density function $f_{\beta, a_1}$}\lbl{bipro}
According to the famous characterization of Dumitriu and Edelman in \cite{DE02}, we know that $(\mu_i)_{1\le i\le m}$ could be regarded as the eigenvalues of the matrix 
${\bf AA}^{\prime},$ where ${\bf A}$ is given as 
\begin{align}\label{matrixa}
\f A=
    \begin{pmatrix}
      x_{1}&  &  &   &  \\
      y_{2}& x_{2} &   &   &\\
       &  y_{3}& x_{3} &   &\\
      &    &   \ddots &\ddots & \\
      &   &   &   y_{m} & x_{m}
    \end{pmatrix}
\end{align} 
with the non-negative random variables $\{x_i\}_{1\le i\le m}$ and $\{y_i\}_{2\le i\le m}$ obeying the distribution and relationships
\begin{itemize}
\item[1).] $\{x_1, x_2, \cdots, x_m, y_2, y_3, \cdots, y_{m}\}$ mutually independent;
\item[2).] $x_i^2\sim \chi^2_{(2a_1-\beta (i-1))}\quad \text{and}\quad  y_i^2\sim\chi^2_{\beta(m-(i-1))}.$
\end{itemize}

Since our calculus below will heavily depend on the properties of $\chi^2$-distribution,  we present Lemma 2.8 in \cite{JM2017}.
\blem\lbl{chi} Given a random variable $X\sim\chi^2_{n}$ for any $n\ge 1.$ Then we have 
$$\aligned & \mathbb{E}X^k=\prod_{l=0}^{k-1}(n+2l), \quad \forall  k\ge 1;  \\
&  \mathbb{E}(X-n)^2=2n; \\
&\mathbb{E}(X-n)^3=8n; \\
&\mathbb{E}(X-n)^4=12n(n+4);\\
&{\rm Var}(X^2)=8n(n+2)(n+3); \\
& {\rm Var}((X-n)^2)=8n(n+6).
\endaligned $$
\nlem 

Now we present two key Lemmas, whose proof are relatively long and will be postponed to the appendix. 
\blem\lbl{key-OO}  Let $\mu=(\mu_1, \mu_2, \cdots, \mu_m)$ be the random variables having joint distribution density $f_{\beta, a_1}$ given in \eqref{bldensity}. We have 
\be\lbl{Var}\aligned  &{\rm Var}\big(\sum_{i=1}^m \mu_i\big)=4a_1m;\\
&\mathbb{E}\sum_{i=1}^m (\mu_{i}-2a_1)^2=4a_1mr;\\
&{\rm Var}\big(\sum_{j=1}^m(\mu_i-2a_1)^2\big)=16\beta a_1m(m-1)(a_1+5)+8 \beta^2 a_1m(m-1)(2m-3)\\
&\quad\quad\quad \quad\quad\quad \quad\quad\quad \quad
+32a_1m(a_1+3);\\
&{\rm Cov}\big(\sum_{i=1}^m(\mu_i-2a_1), \sum_{i=1}^m(\mu_i-2a_1)^2\big)=16a_1mr;\\
&\mathbb{E}\sum_{i=1}^m(\mu_i-2a_1)^3=2\beta^2a_1m(m-1)(m-2)
+12\beta a_1 m(m-1)+16a_1m
\endaligned \ee
for $m\ge 2$ and $a_1>\frac{\beta }{2}(m-1).$  
\nlem 

\blem\lbl{key-O-3}  
Let $\mu=(\mu_1, \mu_2, \cdots, \mu_m)$ be the random variables having joint distribution density $f_{\beta, a_1}$ given in \eqref{bldensity}. Suppose 
that $a_1=O(m),$ then we have 
$$  
{\rm Var}\big(\sum_{i=1}^m(\mu_i-2a_1)^3\big)=O(m^7)
$$ 
for $m$ sufficiently large. 
\nlem 

Next we give two lemmas to describe the property of $\max_{1\le i\le m}|\mu_i-2a_1|$ under the assumption  {\bf A2} or {\bf A3}, respectively. 

\begin{lem}\lbl{limit0o} Let $\mu=(\mu_1, \mu_2, \cdots, \mu_m)$ be the random variables having joint distribution density $f_{\beta, a_1}$ given in \eqref{bldensity}. Suppose $a_1$ and $m$ satisfy $m\to\infty$ and $\frac{m}{a_1}\to 0,$  then  $\max_{1\leq i \leq m}|\frac{\mu_i}{2a_1}-1|\stackrel{p}{\to} 0$  as $m\to\infty$.
\end{lem}
\bprf[\bf Proof.]  Review (1.2) from \cite{JiangLi15}. Treat  $n$ as our ``$m$"  in Theorems 2 and 3 from \cite{JiangLi15}. The rate function $I$ satisfies $I(1)=0$ in both Theorems. 
By  the large deviations in the two Theorems, we see
\beaa
\frac{1}{2a_1}\max_{1\leq i \leq m}\mu_i \stackrel{p}{\to} 1\ \ \mbox{and}\ \ \frac{1}{2a_1}\min_{1\leq i \leq m}\mu_i \stackrel{p}{\to} 1
\eeaa
as $a_1\to\infty$. The conclusion then follows from the inequality
\beaa
\max_{1\leq i \leq m}|\frac{\mu_i}{2a_1}-1| \leq \Big|\frac{\max_{1\leq i \leq m}\mu_i}{2a_1}-1\Big| + \Big|\frac{\min_{1\leq i \leq m}\mu_i}{2a_1}-1\Big|.
\eeaa
The proof is complete. \nprf 

\blem\lbl{maxabso} Let $(\mu_1, \mu_2, \cdots, \mu_m)$ be as in the setting of Lemma \ref{limit0o}. Suppose the assumption {\bf A3} holds, i.e., 
$$ \lim_{a_2\to\infty}\frac{a_1}{\sqrt{a_2}}=x>0 \quad \text{and} \quad \lim_{a_2\to\infty}\frac{m}{\sqrt{a_2}}=y>0.$$
 Then we have 
$$\max_{1\le i\le m}|\frac{\mu_i-2a_1}{m}|\le (1+2\gamma^{-1/2})\beta+1$$ with probability one and $$\max_{1\le i\le m}|\frac{\mu_i-2a_1}{2a_2}|\stackrel{p}{\to} 0$$
as $a_2\to\infty$ with $\gamma:=\frac{\beta y}{2x}.$
\nlem  
\bprf[\bf Proof.] Set $\mu_{\rm max}=\max_{1\leq i \leq m}\mu_i$ and $\mu_{\rm min}=\min_{1\leq i \leq m}\mu_i.$
We get $\gamma\in (0, 1]$  from the condition $2a_1>\beta (m-1).$  Therefore Theorem 10.2.2 in \cite{Dthesis} tells 
$$\frac{\mu_{\rm max}}{m}\to\beta (1+\sqrt{\gamma^{-1}})^2 \quad\text{and}\quad  \frac{\mu_{\rm min}}{m}\to\beta (1-\sqrt{\gamma^{-1}})^2 \quad {\rm a.s.}$$
as $m\to\infty.$
It entails that 
$$ (1+2\gamma^{-1/2})\beta-\frac{1}{2}\le \frac{\mu_{\rm min}}{m}-\frac{\beta}{\gamma}\le \frac{\mu_i}{m}-\frac{\beta}{\gamma}\le \frac{\mu_{\rm max}}{m}-\frac{\beta}{\gamma}\le (1+2\gamma^{-1/2})\beta+\frac{1}{2},$$
almost surely for $1\le i\le m.$ This implies with probability one that as $m$ large enough
$$\max_{1\le i\le m}|\frac{\mu_i}{m}-\frac{\beta}{\gamma}|\le \max\bigg\{|\frac{\mu_{\rm max}}{m}-\frac{\beta}{\gamma}|, \; |\frac{\mu_{\rm min}}{m}-\frac{\beta}{\gamma}|\bigg\}\le \frac{1}{2}+(1+2\gamma^{-1/2})\beta.$$
It is trivial that 
$$|\frac{\mu_i-2a_1}{m}|\le |\frac{\mu_i}{m}-\frac{\beta}{\gamma}|+|\frac{2a_1}{m}-\frac{\beta}{\gamma}|
$$ 
for $1\le i\le m.$
Obviously, by the assumption  {\bf A3},  it holds $|\frac{2a_1}{m}-\frac{\beta}{\gamma}|\le 
\frac12$ as $m\to\infty.$  Immediately one gets with probability one
$$\max_{1\le i\le m}|\frac{\mu_i-2a_1}{m}|\le (1+2\gamma^{-1/2})\beta+1$$
as $m\to\infty$. The proof is then complete since
$$\max_{1\le i\le m}|\frac{2a_1-\mu_i}{2a_2}|=\frac{m}{2a_2}\max_{1\le i\le m}|\frac{\mu_i-2a_1}{m}|\le \big((1+2\gamma^{-1/2})\beta+1\big)\frac{m}{a_2}.$$  
\nprf 

\section{Proof of Proposition \ref{uniqprop}}\lbl{prfprop} 
\subsection{Proof of Proposition \ref{uniqprop} under {\bf A2}.}\lbl{centralA2}  

Review 
$$U_m:=\frac{r}{2a_2}\sum_{i=1}^m (\mu_i-2a_1)-\frac{(a_2-r)}{8a_2^2}\sum_{i=1}^m(\mu_i-2a_1)^2$$
and  the assumption  {\bf A2} : 
$$ m\to\infty, \quad \lim_{a_2\to\infty}\frac{m}{a_1}=0 \quad \text{and}\quad \lim_{a_2\to\infty}\frac{a_1m}{a_2}=\sigma.  
$$
By Lemma \ref{key-OO} and \eqref{cov}, we know $\mathbb{E}\sum_{i=1}^m(\mu_i-2a_1)=0$ and 
$$\frac{r^2}{a_2^2}{\rm Var}(\sum_{i=1}^m\mu_i)=\frac{4r^2a_1m}{a_2^2}=\frac{(2+\beta (m-1))^2a_1m}{a_2^2}\to 0$$
as $a_2\to\infty,$ which is guaranteed by the assumption {\bf A2}.  
This means 
$$\frac{r}{2a_2}\sum_{i=1}^m (\mu_i-2a_1)\stackrel{p}{\to} 0$$
as $a_2\to\infty.$
It remains to prove that 
\be\lbl{secondton}\frac{(a_2-r)}{8a_2^2}\sum_{i=1}^m\big((\mu_i-2a_1)^2-4a_1mr\big)\to N(0, \frac{\beta\sigma^2}{4})\ee
weakly as $a_2\to\infty.$ 
Review the expression \eqref{sumofsquare}: 
$$\sum_{i=1}^m(\mu_i-2a_1)^2=\sum_{i=1}^{m}(z_i-2a_1)^2+2\sum_{i=2}^{m}x_{i-1}^2y_i^2,
$$
where $z_i:=x_i^2+y_i^2$ for $1\le i\le m$ with $y_1:=0.$ 
By \eqref{var-expression} in the appendix, we know  
$${\rm Var}\bigg(\frac{a_2-r}{a_2^2}\sum_{i=1}^{m}(z_i-2a_1)^2\bigg)=\frac{(a_2-r)^2}{a_2^4}(\frac{16}{3}\beta^2a_1m^3+O(a_1^2m+a_1m^2))\to 0$$
since $\frac{a_1m^3}{a_2^2}=O(\frac{m}{a_1})\to 0$ as $a_2\to\infty.$ 
This claims 
\be\lbl{sec}\frac{(a_2-r)}{8a_2^2}\sum_{i=1}^{m}\bigg((z_i-2a_1)^2-\mathbb{E}(z_i-2a_1)^2\bigg)\stackrel{p}{\to}0\ee
as $a_2\to\infty.$
Also the expression \eqref{inter-expression} leads
$$ {\rm Var}\bigg( \frac{2}{a_1m}\sum_{i=2}^{m}x_{i-1}^2y_i^2\bigg)=\frac{4}{a_1^2m^2}\big( 4\beta a_1^2 m^2+o(a_1^2m^2)\big)=16\beta+o(1).
$$
Moreover, the decomposition $$x_{i-1}^2y_{i}^2-\mathbb{E}x_{i-1}^2y_{i}^2=y_i^2(x_{i-1}^2-\mathbb{E}x_{i-1}^2)+(y_{i}^2-\mathbb{E}y_{i}^2)\mathbb{E}x_{i-1}^2$$ provides
$$\aligned \mathbb{E}(x_{i-1}^2y_{i}^2-\mathbb{E}x_{i-1}^2y_{i}^2)^4&\le 8\mathbb{E}y_{i}^8\mathbb{E}(x_{i-1}^2-\mathbb{E}x_{i-1}^2)^4+8(\mathbb{E}x_{i-1}^2)^4\mathbb{E}(y_{i}^2-\mathbb{E}y_{i}^2)^4. \\
\endaligned $$
It follows from Lemma \ref{chi} that the right hand side has the same order as 
$$ (\mathbb{E}y_i^2)^4(\mathbb{E}x_{i-1}^2)^2+(\mathbb{E}x_{i-1}^2)^4(\mathbb{E}y_i^2)^2=O(a_1^4m^2+a_1^2m^4).$$
This means 
\be\lbl{zorder}\frac{1}{a_1^4m^4}\sum_{i=2}^m\mathbb{E}(x_{i-1}^2y_{i}^2-\mathbb{E}x_{i-1}^2y_{i}^2)^4=\frac{O(a_1^4m^3)}{a_1^4m^4}=O(\frac{1}{m})\to 0\ee
as $m\to\infty.$  
Thus by Lyapunov central limit theorem for the sum of independent random variables, we know 
\be\lbl{first}\frac{2}{a_1m}\sum_{i=2}^{m}(x_{i-1}^2y_i^2-\mathbb{E}x_{i-1}^2y_i^2)\to N(0, 16\beta)\ee
weakly as $a_2\to\infty.$
Therefore \eqref{sec} and \eqref{first} and the condition $\lim_{a_2\to\infty}\frac{a_1m}{a_2}=\sigma$ establish \eqref{secondton}. 
The proof of Proposition \ref{uniqprop} under {\bf A2} is complete now. 
 
\hfill$\square$
\subsection{Proof of Proposition \ref{uniqprop} under {\bf A3}.}\lbl{centralA3} 
Recall the assumption {\bf A3}: 
$$ \lim_{a_2\to\infty}\frac{m}{\sqrt{a_2}}=y \quad\text{and}\quad \lim_{a_2\to\infty}\frac{m}{\sqrt{a_2}}=x.
$$
By \eqref{form-single} and \eqref{sumofsquare}, $U_m+\frac{(a_2-r)a_1mr}{2a_2^2}$ could be rewritten as   
$$\aligned U_m+\frac{(a_2-r)a_1mr}{2a_2^2}
&=\frac{r}{2a_2}\sum_{i=1}^m(z_i-\mathbb{E}z_i)-\frac{a_2-r}{4a_2^2}\sum_{i=2}^{m}(x_{i-1}^2y_{i}^2-\mathbb{E}x_{i-1}^2y_{i}^2)\\
&-\frac{(a_2-r)}{8a_2^2}\sum_{i=1}^m\bigg((z_i-2a_1)^2-\mathbb{E}(z_i-2a_1)^2\bigg). \endaligned $$
Set $X_i=Y_i-Z_i,$ where $Y_{i}:=Y_{1, i}-Y_{2, i}$ with   
$$Y_{1, i}:=\frac{r}{2a_2}(z_i-\mathbb{E}z_i), \quad Y_{2, i}:=\frac{(a_2-r)}{8a_2^2}\bigg((z_i-2a_1)^2-\mathbb{E}(z_i-2a_1)^2\bigg)$$
and $$Z_i:=\frac{a_2-r}{4a_2^2}(x_{i-1}^2y_{i}^2-\mathbb{E}x_{i-1}^2y_{i}^2).$$
By the independence of $\{x_i\}_{1\le i\le m}$ and $\{y_j\}_{2\le j\le m},$ we know both $(Y_i)_{2\le i\le m}$ and $(Z_{i})_{2\le i\le m}$ are independent sequences and moreover $Y_i$ is independent of 
$Z_j$ once $j\neq i$ and $j\neq i+1.$
This ensures that 
$(X_i)_{1\le i\le m}$ is a $1$-dependent random variable sequences. Precisely, $X_i$ is independent of $X_j$ for any $j$ satisfying $|j-i|>1.$ 

Now we follow the idea in \cite{HR} to separate the sum into two parts, both of which are the sum of independent random variables. We will prove that one of them tends to zero in probability and the other one tends to a normal distribution weakly as $a_2\to\infty.$ 
For that aim, we choose $\kappa=[m^{\alpha}]$ with  $0<\alpha<1$ and $\nu=[\frac{m}{\kappa}].$ Then $m=\kappa \nu+s$ with  $0\le s<\kappa.$ Obviously $\frac{\nu}{m}\to 0$ and  $\frac{\kappa}{m}\to 0$ as $m\to\infty.$    
For any $1\le i\le \nu,$ set
$$ W_i:=\sum_{l=1}^{\kappa-1}X_{(i-1)\kappa+l} \quad\text{and}\quad V_i:=X_{i\kappa} 
$$
and $W_{\nu+1}:=\sum_{l=1}^{s} X_{\nu\kappa+l}.$
Then both $(W_i)_{1\le i\le \nu+1}$ and $(V_i)_{1\le i\le \nu}$ are independent random variable sequences and 
$$\sum_{i=1}^m X_i=\sum_{i=1}^{\nu+1} W_i+\sum_{i=1}^{\nu}V_{i}.$$ Next we will prove that 
\be\lbl{key-centr} 
\sum_{i=1}^{\nu}V_{i}\stackrel{p}{\to} 0 \quad\text{and} \quad \sum_{i=1}^{\nu+1} W_i\to N\big(0, \;  \frac{\beta \sigma^2}{4}\big)
\ee
weakly as $a_2\to\infty.$ 
Once \eqref{key-centr} holds, the proof of Proposition \ref{uniqprop} is complete. By definition and the property of variance, it follows  
$${\rm Var}(V_{i})\le 3{\rm Var}(Y_{1, i\kappa})+3{\rm Var}(Y_{2, i\kappa})+3{\rm Var}(Z_{i\kappa}).$$
Easily it holds
$${\rm Var}(Y_{1, i\kappa})=\frac{r^2}{4a_2^2}{\rm Var}(z_{i\kappa})=\frac{r^2}{2a_2^2}(2a_1+b_{i\kappa}).$$
Based on \eqref{varz} and \eqref{varinter}, we have 
$${\rm Var}(Y_{2, i\kappa})=\frac{(a_2-r)^2}{8a_2^4}\big((2a_1+b_{i\kappa})^2+(2a_1+b_{i\kappa})(b_{i\kappa}^2+2b_{i\kappa}+6)\big)$$
and 
$${\rm Var}(Z_{i\kappa})=\frac{\beta(a_2-r)^2}{8a_2^4}(2a_1-\beta(i\kappa-2))(m+1-i\kappa)(2a_1+b_{i\kappa}+\beta+2).$$
By the assumption {\bf A3}, $m$ and $a_1$ have the same order as $\sqrt{a_2},$ then all these three terms above have order $m^{-1}.$
It follows immediately 
$${\rm Var}(\sum_{i=1}^{\nu}V_{i})= \sum_{i=1}^{\nu} {\rm Var}(V_i)=\nu O(m^{-1}) \to 0$$ 
as $m\to\infty.$  Since $\sum_{i=1}^{\nu}\mathbb{E}V_i=0,$ the first limit in \eqref{key-centr} is verified. 
Now we work harder on the tough second term. Since $\sum_{i=1}^{\nu+1}W_i$ is the sum of independent random variables and $\mathbb{E}(W_i)=0,$ by Lyapunov central limit Theorem again,  it suffices to prove that 
\be\lbl{var-W} \sum_{i=1}^{\nu+1} \mathbb{E} W_i^2\to\frac{\beta \sigma^2}{4} \quad \text{and} \quad 
\sum_{i=1}^{\nu+1} \mathbb{E}W_i^4\to 0.\ee 
Since $(X_i)_{1\le i\le m}$ is $1$-dependent and $\mathbb{E}X_i=0,$ we have 
 \be\lbl{square-exp-for}\mathbb{E}(\sum_{i=p}^q X_i)^2=\sum_{i=p}^q \mathbb{E}X_i^2+2\sum_{i=p}^{q-1} \mathbb{E}[X_iX_{i+1}]\ee
  for any integers  
$1\le p<q.$   
 Thereby it follows 
 $$ \aligned  \sum_{i=1}^{\nu+1} \mathbb{E} W_i^2&=\sum_{i=1}^{\nu}\sum_{l=1}^{\kappa-1}\mathbb{E}X_{(i-1)\kappa+l}^2+2\sum_{i=1}^{\nu}\sum_{l=1}^{\kappa-2}\mathbb{E}[X_{(i-1)\kappa+l}X_{(i-1)\kappa+l+1}]\\
&+\sum_{l=1}^{s}\mathbb{E}X_{\nu \kappa+l}^2+2\sum_{l=1}^{s-1}\mathbb{E}X_{\nu\kappa+l}X_{\nu\kappa+l+1}\\
 &=\sum_{i=1}^m\mathbb{E}X_i^2+2\sum_{i=1}^{m-1}\mathbb{E}[X_iX_{i+1}]-\sum_{i=1}^{\nu}\mathbb{E}X_{i\kappa}^2-2\sum_{i=1}^{\nu}\mathbb{E}[(X_{i\kappa-1}+X_{i\kappa+1})X_{i\kappa}].
 \endaligned $$
 As for $${\rm Var}(\sum_{i=1}^{\nu} V_i)={\rm Var}(\sum_{i=1}^{\nu} X_{i\kappa})=\sum_{i=1}^{\nu}\mathbb{E}X_{i\kappa}^2\to 0,$$ we could similarly have as $m\to\infty$
 $$\sum_{i=1}^{\nu}\mathbb{E}|(X_{i\kappa-1}+X_{i\kappa+1})X_{i\kappa}|\le 2\sum_{i=1}^{\nu}\mathbb{E}X_{i\kappa}^2+\sum_{i=1}^{\nu}\mathbb{E}X_{i\kappa-1}^2+\sum_{i=1}^{\nu}\mathbb{E}X_{i\kappa+1}^2\to 0.$$
 Therefore for  the first limit in \eqref{var-W}, it remains to prove  
 \be\lbl{last} \sum_{i=1}^m\mathbb{E}X_i^2+2\sum_{i=1}^{m-1}\mathbb{E}[X_iX_{i+1}]=\mathbb{E}(\sum_{i=1}^m X_i)^2\to \frac{\beta \sigma^2}{4}\ee
 as $a_2\to\infty.$
By definition, 
$$ \mathbb{E}(\sum_{i=1}^mX_i)^2=\mathbb{E} \bigg( \frac{r}{2a_2}\sum_{i=1}^m(\mu_i-2a_1)-\frac{a_2-r}{8a_2^2}\sum_{i=1}^m\big[(\mu_i-2a_1)^2-\mathbb{E}(\mu_i-2a_1)^2\big]\bigg)^2.
$$
By Lemma \ref{key-OO} and condition {\bf A3}, certainly we have 
$$\aligned \mathbb{E}(\sum_{i=1}^mX_i)^2&=\frac{r^2}{4a_2^2}\cdot 4a_1m+\frac{(a_2-r)^2}{64a_2^4}\cdot (16\beta a_1^2m^2+16\beta^2 a_1m^3)\\
&-\frac{r(a_2-r)}{8a_2^3}8\beta a_1m^2+o(1)\\
&=\frac{\beta^2xy^3}{4}+\frac{\beta \sigma^2}{4}+\frac{\beta^2xy^3}{4}-\frac{\beta^2xy^3}{2}+o(1)\\
&=\frac{\beta \sigma^2}{4}+o(1)
\endaligned $$ 
as $a_2$ large enough. 
Therefore \eqref{last} is satisfied.  
The last thing left is to verify the second limit in \eqref{var-W}.  

For any $1\le p<q\le m,$
we have 
$$\aligned \mathbb{E}(\sum_{i=p}^q X_i)^4&=\sum_{i=p}^{q}\mathbb{E}X_i^4+3\sum_{p\le i\neq j\le q}\mathbb{E}X_{i}^2X_{j}^2+4\sum_{p\le i\neq j\le q}\mathbb{E}X_{i}^3X_{j}\\
&+6\sum_{p\le i\neq j\neq k\le q}\mathbb{E}(X_i^2X_jX_k)+\sum_{p\le i\neq j\neq k\neq l\le q}\mathbb{E}X_iX_jX_kX_l.\\
\endaligned $$
Since $(X_i)_{1\le i\le m}$ is $1$-dependent and $\mathbb{E}X_i=0$ for $1\le i\le m,$ we know 
\be\lbl{final-up}\aligned \mathbb{E}(\sum_{i=p}^q X_i)^4&=\sum_{i=p}^{q}\mathbb{E}X_i^4+6\sum_{p\le i\le j-1\le q-1}\mathbb{E}[X_{i}^2X_{j}^2]+4\sum_{i=p}^{q-1}\mathbb{E}[X_{i}X_{i+1}(X_i^2+X_{i+1}^2)]\\
&\le \sum_{i=p}^{q}\mathbb{E}X_i^4+6\sum_{p\le i\le j-1\le q-1}(\mathbb{E}X_{i}^4)^{1/2}(\mathbb{E}X_{j}^4)^{1/2}\\
&+4\sum_{i=p}^{q-1}(\mathbb{E}X_{i}^4)^{1/4}(\mathbb{E}X_{i+1}^4)^{3/4}+4\sum_{i=p}^{q-1}(\mathbb{E}X_{i}^4)^{3/4}(\mathbb{E}X_{i+1}^4)^{1/4}.
\endaligned \ee
Now we investigate the dominated order of $\mathbb{E}X_i^4$ for $1\le i\le m.$ 
With the decomposition $X_i=Y_{1, i}-Y_{2, i}-Z_{i},$ we know 
\be\lbl{final-mix}\mathbb{E}X_i^4\le 27\mathbb{E}Y_{1, i}^4+27\mathbb{E}Y_{2, i}^4+27\mathbb{E}Z_{i}^4.\ee
As in \eqref{zorder},  it holds 
\be\lbl{final-1}\mathbb{E}Z_{i}^4=\frac{1}{a_2^4}O(a_1^4m^2+a_1^2m^4)=O(\frac{1}{m^2}).\ee
Since $z_{i}\sim\chi^2_{2a_1+b_i}$ with $b_i=\beta (m-2(i-1)),$ we have from Lemma \ref{chi}
\be\lbl{final-2}\mathbb{E}Y_{1, i}^4=\frac{r^4}{(2a_2)^4}\mathbb{E}(z_i-\mathbb{E}z_i)^4=O(\frac{r^4(2a_1+b_i)^2}{a_2^4})=O(\frac1{m^2}).\ee
Now we work on the term $\mathbb{E}Y_{2, i}^4.$ Indeed 
$$\aligned &\mathbb{E}\big((z_i-2a_1)^2-\mathbb{E}(z_i-2a_1)^2\big)^4\\
=&\mathbb{E}\big((z_i-\mathbb{E}z_i)^2-{\rm Var}(z_i)+2b_i(z_i-\mathbb{E}z_i)\big)^4\\
\le &27 \mathbb{E}(z_i-\mathbb{E}z_i)^8+27({\rm Var}(z_i))^4+432 b_i^4\mathbb{E}(z_i-\mathbb{E}z_i)^4.\\
\endaligned $$
The second term has order $m^4$ and the third part has order $m^6$ from the property of chi square distribution and the condition $a_1=O(m).$ 
Suppose $X\sim\chi^2_{n}.$ By binomial expansion, 
$$\aligned (X-n)^8&=X^8-8nX^7+28n^2X^6-56n^3X^5+70n^4X^4
-56n^5X^3+28n^6X^2-8n^7X+n^8.\\
\endaligned $$
Applying Lemma \ref{chi} to $\mathbb{E}X^k$ for $1\le k\le 8,$ and combining carefully alike terms, we finally have 
$$\mathbb{E}(X-n)^8=7\cdot 240n^4+o(n^4).$$  This means 
$\mathbb{E}(z_i-\mathbb{E}z_i)^8=O((2a_1+b_i)^4)=O(m^4).$ 
Therefore 
\be\lbl{final-3}\mathbb{E}Y_{2, i}^4=\frac{O(m^6)}{a_2^4}=O(\frac1{m^2}).\ee
Putting \eqref{final-1}, \eqref{final-2} and \eqref{final-3} back into \eqref{final-mix}, we know that 
$$\mathbb{E}X_i^4=O(\frac1{m^2})$$
for all $1\le i\le m.$ 
This and \eqref{final-up} tell that 
$$ \mathbb{E}(\sum_{i=p}^q X_i)^4=O(\frac{1}{m^2})(q-p)^2.
$$ 
Then 
$$\sum_{i=1}^{\nu}\mathbb{E}W_i^4=\sum_{i=1}^{\nu}\mathbb{E}(\sum_{l=1}^{\kappa-1}X_{(i-1)\kappa+l})^4=\nu \kappa^2 O(\frac{1}{m^2})=O(\frac{\kappa}{m})\to 0 $$
as $m\to\infty.$ This is exactly the second limit in \eqref{var-W}. The proof is complete now. 
\hfill$\square$
\bcor\label{A3partial}  Let $\mu=(\mu_1, \mu_2, \cdots, \mu_m)$ be random variables with density $f_{\beta, a_1}$ as in \eqref{bldensity}. 
Then under the assumption {\bf A3}, with $\sigma:=xy$ in {\bf A3} we have 
$$\frac{(a_2-r)}{8a_2^2}\sum_{i=1}^m (\mu_i-2a_1)^2\to N(0, \frac{\beta \sigma^2}{4}+\frac{\beta^2 xy^3}{4})$$
weakly as $a_2\to\infty.$
\ncor 
\bprf  Let $\nu, \kappa, Y_{2, i}$ and $Z_i$ be the same as in the proof of Proposition \ref{uniqprop}. Set $\tilde{X}_i=Y_{2, i}+Z_i$ for $1\le i\le m$  and  $$\tilde{W}_i=\sum_{l=1}^{\kappa-1}\tilde{X}_{(i-1)\kappa+l} \quad  \text{and} \quad \tilde{V}_i=\tilde{X}_{i\kappa}$$ for any $1\le i\le \nu$ and $\tilde{W}_{\nu+1}=\sum_{l=1}^s \tilde{X}_{\nu\kappa+l}.$ Then by \eqref{sumofsquare}, we know  $$\frac{a_2-r}{8a_2^2}\sum_{i=1}^m( (\mu_i-2a_1)^2-\mathbb{E}(\mu_i-2a_1)^2)=\sum_{i=1}^m\tilde{X}_i=\sum_{i=1}^{\nu}\tilde{V}_i+\sum_{i=1}^{\nu+1}\tilde{W_i}.$$
According to the argument above, we could prove that 
$$\sum_{i=1}^{\nu} \tilde{V}_i\stackrel{p}{\to} 0 \quad\text{and}\quad \mathbb{E}\sum_{i=1}^{\nu+1}\tilde{W}_i^4\to 0$$ 
as $a_2\to\infty.$ The only difference is 
$$\aligned \mathbb{E}(\sum_{i=1}^{m} \tilde{X}_i)^2&=\frac{(a_2-r)^2}{64a_2^4} {\rm Var}\big(\sum_{i=1}^m (\mu_i-2a_1)^2\big)\\
&=\frac{(a_2-r)^2}{64a_2^4}\big(16\beta a_1^2m^2+16\beta^2 a_1m^3+o(a_1^2m^2+a_1m^3)\big),
\endaligned $$ 
where for the second equality we use Lemma \ref{key-OO}. Therefore  by the assumption {\bf A3},  it follows $$\mathbb{E}(\sum_{i=1}^{m} \tilde{X}_i)^2\to \frac{\beta \sigma^2+\beta^2xy^3}{4}$$
as $a_2\to\infty.$ 
The same limit holds for $\sum_{i=1}^{\nu+1}\mathbb{E}\tilde{W}_i^2.$ Since $\mathbb{E}\tilde{W}_i=0$ for $1\le i\le \nu+1,$ finally we know 
$$\sum_{i=1}^{\nu+1}\tilde{W}_i\to N(0,  \frac{\beta \sigma^2+\beta^2xy^3}{4})$$
weakly as $a_2\to\infty.$ Hence 
$$\frac{a_2-r}{8a_2^2}\sum_{i=1}^m( (\mu_i-2a_1)^2-\mathbb{E}(\mu_i-2a_1)^2)=\sum_{i=1}^{\nu}\tilde{V}_i+\sum_{i=1}^{\nu+1}\tilde{W_i}{\to} N(0,  \frac{\beta \sigma^2+\beta^2xy^3}{4})$$ weakly as $a_2\to\infty.$ The proof is complete. 
\nprf

\section{Proof of Theorem \ref{bjconv}}\lbl{prfthm} 
In this section, we will give the final statement on the proof of Theorem \ref{bjconv}. 

Recall the joint density function of $\lambda=(\lambda_1, \lambda_2, \cdots, \lambda_m)\in [0, 1]^m$ is given by 
\be\lbl{bjdensityf}
f_{\beta, a_1, a_2}(x_1, x_2, \cdots, x_m)=C_{\rm J}^{\beta, a_1, a_2}\prod_{1\le i<j\le m}|x_i-x_j|^{\beta}\prod_{i=1}^m x_i^{a_1-r}(1-x_i)^{a_2-r},
\ee
where $a_1, a_2>\beta (m-1)/2$ and $r=1+\frac{\beta}{2}(m-1),$ and
$$C_{\rm J}^{\beta, a_1, a_2}=\prod_{j=1}^m \dfrac{\Gamma(1+\beta/2)\Gamma(a_1+a_2-(\beta/2)(m-j))}{\Gamma(1+(\beta/2)j)\Gamma(a_1-(\beta/2)(m-j))\Gamma(a_2-(\beta/2)(m-j))}.$$
It is clear that the joint distribution density for $\theta:=2a\lambda$, denoted by $g_{\beta, a_1, a_2},$  should be as follows
$$\aligned & g_{\beta, a_1, a_2}(x_1, x_2, \cdots, x_m)\\
&:=C_{\rm J}^{\beta, a_1, a_2}(\frac{1}{2a})^{c}\prod_{1\le i<j\le m}|x_i-x_j|^{\beta}\prod_{i=1}^m x_i^{a_1-r}(1-\frac{x_i}{2a})^{a_2-r}{\bf I}_{\{\max\theta_i\le 2a\}}\\
&=C_{\rm J}^{\beta, a_1, a_2}(\frac{1}{2a})^{a_1m}\prod_{1\le i<j\le m}|x_i-x_j|^{\beta}\prod_{i=1}^m x_i^{a_1-r}(1-\frac{x_i}{2a})^{a_2-r} {\bf I}_{\{\max\theta_i\le 2a\}} ,
\endaligned $$ 
where $c:=\beta m(m-1)/2+m(a_1-r)+m=a_1m.$ 
Review the joint density function of $\mu=(\mu_1, \mu_2, \cdots, \mu_m)$ is
\be\lbl{bldensityf}
f_{\beta, a_1}(x_1, x_2, \cdots, x_m)=C_{\rm L}^{\beta, a_1}\prod_{1\le i<j\le m}|x_i-x_j|^{\beta}\prod_{i=1}^m x_i^{a_1-r}e^{-1/2\sum_{i=1}^m x_i},
\ee
where $$C_{\rm L}^{\beta, a_1}=2^{-ma_1}\prod_{j=1}^m \dfrac{\Gamma(1+\beta/2)}{\Gamma(1+(\beta/2)j)\Gamma(a_1-(\beta/2)(m-j))}.$$
Remember 
\be\lbl{KLf}\aligned
& K_m=(\frac{1}{a})^{m a_1}\prod_{i=0}^{m-1} \frac{\Gamma(a-\eta i)}{\Gamma(a_2-i\eta)}; \\
& L_m(\mu)=e^{\frac12\sum_{i=1}^m\mu_i}\prod_{i=1}^{m}(1-\frac{\mu_i}{2a})^{a_2-r} {\f I}_{\{\max \mu_i\le 2a\}}.
\endaligned \ee
Observing the expressions \eqref{bjdensityf}, \eqref{bldensityf} and \eqref{KLf}, we know 
$$\frac{g_{\beta, a_1, a_2}}{f_{\beta, a_1}}=K_mL_m.$$
This leads 
$$\aligned \|\ml{L}(2a\lambda)-\ml{L}(\mu)\|_{\rm TV}&=\int_{[0, +\infty)^m}|g_{\beta, a_1, a_2}(x)-f_{\beta, a_1}(x)|dx\\
&=\int_{[0, \infty)^m}|\frac{g_{\beta, a_1, a_2}(x)}{f_{\beta, a_1}(x)}-1| f_{\beta, a_1}(x) dx\\
&=\mathbb{E}|K_mL_m(\mu)-1|. \endaligned 
$$ 
Meanwhile, the Kullback-Leibler distance $D_{\rm KL}\big(\mathcal{L}(a \lambda)||\mathcal{L}(\mu)\big)$ could be expressed as 
$$\aligned D_{\rm KL}\big(\mathcal{L}(2a \lambda)||\mathcal{L}(\mu)\big)&=\int_{[0, \infty)^m}\frac{g_{\beta, a_1, a_2}(x)}{f_{\beta, a_1}(x}\log \frac{g_{\beta, a_1, a_2}(x)}{f_{\beta, a_1}(x)}f_{\beta, a_1}(x)dx\\
&=\int_{[0, \infty)^m}\log \frac{g_{\beta, a_1, a_2}(x)}{f_{\beta, a_1}(x)}g_{\beta, a_1, a_2}(x)dx\\
&=\mathbb{E}\log(K_mL_m(\lambda)).
\endaligned $$

As in \cite{JM2017}, we consider another modified version $L_m^{\prime}, K_m^{\prime}$ of $L_m, K_m$ respectively,  which are defined by 
 \be\lbl{Lprime}\aligned  &L_m^{\prime}=(1+\frac{a_1}{a_2})^{m(a_2-r)} L_m; \\
 &K_m^{\prime}= (1+\frac{a_1}{a_2})^{-m(a_2-r)}K_m.
\endaligned  \ee 
Obviously $L_m^{\prime}K_m^{\prime}=L_mK_m.$ Therefore we have 
\be\lbl{mexpre}\aligned  &\|\ml{L}(2a\lambda)-\ml{L}(\mu)\|_{\rm TV}=\mathbb{E}|K_m^{\prime}L_m^{\prime}(\mu)-1|;\\
&D_{\rm KL}\big(\mathcal{L}(2a \lambda)||\mathcal{L}(\mu)\big)=\mathbb{E}\log(K_m^{\prime}L_m^{\prime}(\lambda)).
\endaligned
\ee

\subsection{Proof of (i) of Theorem \ref{bjconv}} By the relationship \eqref{TV_KL} mentioned in the introduction, to prove (i) of Theorem \ref{bjconv}, we just need to prove 
$$\lim_{a_2\to\infty}D_{\rm KL}\bigg(\mathcal{L}(2a \lambda)||\mathcal{L}(\mu)\bigg)=0. $$ 
By Lemma \ref{kmO}, since $a_1m=o(a_2)$ and $\beta(m-1)< 2a_1,$ then one gets $a_1m^3=o(a_2^2).$ Recalling $r=\eta (m-1)+1,$ we have 
$$\aligned \log K_m^{\prime}&=\log(K_m)-m(a_2-r)\log(1+\frac{a_1}{a_2})\\
&=-a_1m+\frac{rm}{2}\log(1+\frac{a_1}{a_2})+o(1)\\
&=-a_1m+\frac{\eta a_1m^2}{2a_2}+o(1).\\
\endaligned $$
Meanwhile by \eqref{KLf} and \eqref{Lprime}, we have 
\be\lbl{mLprime}\aligned \mathbb{E}\log(L_m^{\prime}(\lambda))&=m(a_2-r)\log(1+\frac{a_1}{a_2})+\frac12\mathbb{E}\sum_{i=1}^{m} \theta_i+(a_2-r)\mathbb{E}\sum_{i=1}^m \log(1-\frac{\theta_i}{2a})\\
&=\frac12\mathbb{E}\sum_{i=1}^{m} \theta_i+(a_2-r)\sum_{i=1}^m\mathbb{E} \log(1+\frac{2a_1-\theta_i}{2a_2}).
\endaligned \ee
Here we use the fact
$$\log(1+\frac{a_1}{a_2})+\log(1-\frac{\theta_i}{2a})=\log(1+\frac{2a_1-\theta_i}{2a_2}).$$
Therefore it follows from \eqref{mexpre} that
\bea\lbl{K-Lfinal} \aligned &D_{\rm KL}\big(\mathcal{L}(2a \lambda)||\mathcal{L}(\mu)\big)=\mathbb{E}\log K_m^{\prime}+\mathbb{E}\log(L_m^{\prime}(\lambda)) \\
&=-a_1m+\frac{\eta a_1m^2}{2a_2}+\frac12\mathbb{E}\sum_{i=1}^{m} \theta_i
+(a_2-r)\mathbb{E}\sum_{i=1}^m \log(1+\frac{2a_1-\theta_i}{2a_2})+o(1)\\
&\le-a_1m+\frac{\eta a_1m^2}{2a_2}+\frac12\mathbb{E}\sum_{i=1}^{m} \theta_i+o(1)\\
&+(a_2-r)\mathbb{E}\sum_{i=1}^m \big(\frac{2a_1-\theta_i}{2a_2}-\frac{(2a_1-\theta_i)^2}{8a_2^2}+\frac{(2a_1-\theta_i)^3}{24a_2^3}\big),\endaligned
\eea
where the last inequality is due to the elementary inequality 
$$\log(1+x)\le x-\frac{x^2}{2}+\frac{x^3}{3}, \quad x>-1.$$
Obviously, applying Proposition \ref{objsquare} to $\dfrac{\theta}{2a},$  we have  by $a_1m=o(a_2)$ that 
$$\aligned &-a_1m+\frac{1}{2}\mathbb{E}\sum_{i=1}^m\theta_i+\frac{a_2-r}{2a_2}\mathbb{E}\sum_{i=1}^m (2a_1-\theta_i)\\
&=-a_1m+\frac{r}{2a_2}\mathbb{E}\sum_{i=1}^m \theta_i+\frac{(a_2-r)a_1m}{a_2}\\
&=-a_1m+\frac{(a_2-r)a_1m}{a_2}+\frac{r}{a_2}a_1m+o(1)\\
&=o(1) \\
\endaligned $$ 
and also 
$$\aligned &\frac{a_2-r}{8a_2^2}\mathbb{E}\sum_{i=1}^m(2a_1-\theta_i)^2=\frac{a_2-r}{8a_2^2}\bigg(4a_1^2m-4a_1 \mathbb{E}\sum_{i=1}^m \theta_i+\mathbb{E}\sum_{i=1}^m\theta_i^2\bigg)\\
&=\frac{a_2-r}{8a_2^2}(4a_1^2m-4a_1\cdot 2a_1m+4a_1^2m+4\eta a_1m^2)+o(1)\\
&=\frac{\eta a_1m^2}{2a_2}+o(1).
\endaligned $$ 
Similarly we get 
$$\aligned &\frac{a_2-r}{24 a_2^3}\mathbb{E}\sum_{i=1}^m(2a_1-\theta_i)^3=\frac{a_2-r}{24a_2^3}\mathbb{E}\sum_{i=1}^m\bigg(8a_1^3-12a_1^2\theta_i+6a_1\theta_i^2-\theta_i^3\bigg)\\
&=\frac{a_2-r}{24a_2^3}\bigg(8a_1^3 m-12a_1^2 \cdot 2a_1 m+6a_1(4a_1^2 m+4\eta a_1 m^2)\\
&-8a_1^3m-24\eta a_1^2m^2-8\eta^2a_1m^3\bigg)+o(1)\\
&=-\frac{(a_2-r)\eta^2 a_1m^3}{3a_2^3}+o(1)=o(1).
\endaligned $$ 
Therefore plugging all these expressions into \eqref{K-Lfinal}, we have 
$$ \aligned &D_{\rm KL}\big(\mathcal{L}(2a \lambda)||\mathcal{L}(\mu)\big)\le-a_1m+\frac{\eta a_1m^2}{2a_2}+\frac12\mathbb{E}\sum_{i=1}^{m} \theta_i\\
&+(a_2-r)\mathbb{E}\sum_{i=1}^m \big(\frac{2a_1-\theta_i}{2a_2}-\frac{(2a_1-\theta_i)^2}{8a_2^2}+\frac{(2a_1-\theta_i)^3}{24a_2^3}\big) \\
&=\frac{\eta a_1m^2}{2a_2}-\frac{\eta a_1m^2}{2a_2}+o(1)=o(1).\endaligned
$$ 
The desired result is obtained.  $\hfill\square$

\subsection{Proof of (ii) of Theorem \ref{bjconv}} 
We first present a crucial Lemma on central limit theorem for $\log(L_m^{\prime}(\mu))$ with $\mu$ having probability density function $f_{\beta, a_1}$ in \eqref{bldensity}. 
\subsubsection{Central limit theorem for $\log(L_m^{\prime}(\mu))$}
We first present the result when the assumption {\bf A2} is satisfied. 
\blem\lbl{LemmaLprime}  Suppose that 
$(\mu_1, \mu_2, \cdots, \mu_m)$ have joint distribution density \eqref{bldensityf} and let $L_m^{\prime}$ be given as in \eqref{Lprime}.  Under the assumption {\bf A2}, we have  
\be\lbl{lim-prime} \log L_m^{\prime}(\mu)-a_1m +\frac{(a_2-r)a_1mr}{2a_2^2}\to N(0, \frac{\beta \sigma^2}{4})
\ee
weakly as $a_2\to\infty.$
\nlem 
\bprf[\bf Proof]  
By Taylor's formula, there exists some continuous function $h$ such that  
$$\log(1+x)=x-\frac{x^2}{2}+x^3h(x)$$ for all $x>-1.$  Based on 
\eqref{mLprime}, we are able to write  
 \be\lbl{second-proof}\aligned \log L_m'(\mu)&=\frac12\sum_{i=1}^m\mu_i+(a_2-r)\sum_{i=1}^m\bigg(\frac{2a_1-\mu_i}{2a_2}-\frac{(2a_1-\mu_i)^2}{8a_2^2}\bigg)\\
 &\quad +(a_2-r)\sum_{i=1}^m(\frac{2a_1-\mu_i}{2a_2})^3h(\frac{2a_1-\mu_i}{2a_2})\\
 &=a_1m+\frac{r}{2a_2}\sum_{i=1}^m(\mu_i-2a_1)-\frac{a_2-r}{8a_2^2}\sum_{i=1}^m(\frac{2a_1-\mu_i}{2a_2})^2\\
 &\quad +(a_2-r)\sum_{i=1}^m(\frac{2a_1-\mu_i}{2a_2})^3h(\frac{2a_1-\mu_i}{2a_2})\\
 &=a_1m+U_m
+ (a_2-r)\sum_{i=1}^m(\frac{2a_1-\mu_i}{2a_2})^3h(\frac{2a_1-\mu_i}{2a_2}).\\
 \endaligned \ee 
Therefore we have 
$$
\aligned  &\log L_m^{\prime}(\mu)-a_1m+\frac{(a_2-r)a_1mr}{2a_2^2}\\
&=U_m+\frac{(a_2-r)a_1mr}{2a_2^2}+(a_2-r)\sum_{i=1}^m(\frac{2a_1-\mu_i}{2a_2})^3h(\frac{2a_1-\mu_i}{2a_2}).
\endaligned 
$$
 By Proposition \ref{uniqprop}, under the assumption {\bf A2}, it holds 
 $$ U_m+\frac{(a_2-r)a_1mr}{2a_2^2} \to N(0, \frac{\beta \sigma^2}{4})$$
 weakly as $a_2\to\infty.$
Hence to prove \eqref{lim-prime}, it remains to prove 
 \be\lbl{lastprime}
\delta_m:= (a_2-r)\sum_{i=1}^m(\frac{2a_1-\mu_i}{2a_2})^3h(\frac{2a_1-\mu_i}{2a_2}) \to 0
\ee
in probability as $a_2\to\infty.$
The proof of \eqref{lastprime} follows that for (2.17) in \cite{JM2017}. 
Review $\log (1+x)=x-\frac{x^2}{2}+x^3h(x)$ with $h$ being a continuous function on $(-1, \infty)$. Then, $\tau:=\sup_{|x|\leq 1/2}|h(x)|< \infty.$  Hence, by the fact $\frac{a_1}{a_2}\to 0$ from {\bf A2}, we have 
\be\lbl{probab}\aligned 
P(|\delta_m|>\epsilon)&= P\Big(|\delta_m|>\epsilon,\ \max_{1\leq i \leq m}|\frac{2a_1-\mu_i}{2a_2}|\leq \frac{1}{2}\Big) + P\Big(|\delta_m|>\epsilon, \max_{1\leq i \leq m}|\frac{2a_1-\mu_i}{2a_2}|>\frac{1}{2}\Big)\\
& \leq  P\Big(|\delta_m|>\epsilon,\ \max_{1\leq i \leq m}|\frac{2a_1-\mu_i}{2a_2}|\leq \frac{1}{2}\Big) + P\Big(\max_{1\leq i \leq m}|\frac{\mu_i}{2a_1}-1|>\frac{1}{4}\Big) \\
\endaligned \ee
as $a_2$ is sufficiently large.
Under $\max_{1\leq i \leq m}|\frac{2a_1-\mu_i}{2a_2}|\leq \frac{1}{2}$,
$$\aligned 
|\delta_m| &\leq  (2\tau)\cdot \max_{1\leq i \leq m}|\frac{2a_1-\mu_i}{2a_2}|\cdot \frac{a_2-r}{8a_2^2}\sum_{i=1}^m(\mu_i-2a_1)^2\\
&=2\tau \cdot \max_{1\leq i \leq m}|\frac{\mu_i}{2a_1}-1|\cdot\frac{a_1}{a_2} \cdot \frac{a_2-r}{8a_2^2}\big(\sum_{i=1}^m(\mu_i-2a_1)^2-4a_1mr\big)\\
&\quad +2\tau \cdot \max_{1\leq i \leq m}|\frac{\mu_i}{2a_1}-1|\cdot\frac{a_1}{a_2} \cdot \frac{(a_2-r)a_1mr}{2a_2^2},
\endaligned 
$$ 
which tends to $0$ in probability as $a_2\to\infty$ since $\max_{1\leq i \leq m}|\frac{\mu_i}{2a_1}-1|\to 0$ in probability by Lemma \ref{limit0o},   $a_1^2mr(a_2-r)/a_2^3\to\beta\sigma^2/2$ and $$ \frac{a_2-r}{8a_2^2}(\sum_{i=1}^m(\mu_i-2a_1)^2-4a_1mr)\to N(0, \frac{\beta\sigma^2}{4})$$ weakly 
as $a_2\to\infty$ by \eqref{secondton} and the assumption {\bf A2} .  This, \eqref{probab} and Lemma \ref{limit0o} again concludes \eqref{lastprime}.
 \nprf
 
 Now we present the parallel one under the assumption {\bf A3}. 
\blem\lbl{centr-O} Let $(\mu_1, \mu_2, \cdots, \mu_m)$ be the random variables having joint distribution density $f_{\beta, a_1}$ given in \eqref{bldensity} and $L_m'$ be given in \eqref{Lprime}. 
Then under the assumption {\bf A3} with $\sigma=xy$, we have 
\be\lbl{lim-primeO} \log L_m^{\prime}(\mu)-a_1m+\frac{a_1m(a_2-r)r}{2a_2^2}{\to}N(-\frac{\beta^2 xy^3}{12}, \frac{\beta \sigma^2}{4})
\ee
weakly as $a_2\to\infty.$
\nlem 
\bprf[\bf Proof.]  Applying the Taylor formula $\log(1+x)=x-\frac{x^2}{2}+\frac{x^3}{3}+x^4h(x)$ with $h$ a continuous function on 
$(-1, +\infty),$ the same argument as for \eqref{second-proof} leads
$$ \aligned  \log L_m'(\mu)&=a_1m+\frac{r}{2a_2}\sum_{i=1}^m(\mu_i-2a_1)
-\frac{(a_2-r)}{8a_2^2}\sum_{i=1}^m(2a_1-\mu_i)^2\\
&+\frac{(a_2-r)}{24a_2^3}\sum_{i=1}^m(2a_1-\mu_i)^3+(a_2-r)\sum_{i=1}^m(\frac{2a_1-\mu_i}{2a_2})^4h(\frac{2a_1-\mu_i}{2a_2}).\endaligned
$$
Then we have 
 $$ \aligned
&\log L_m'(\mu)-a_1m+\frac{(a_2-r)a_1mr}{2a_2^2}\\
&=U_m+\frac{(a_2-r)a_1mr}{2a_2^2}
-\frac{(a_2-r)}{24a_2^3}\sum_{i=1}^m(\mu_i-2a_1)^3\\
&\quad +(a_2-r)\sum_{i=1}^m(\frac{2a_1-\mu_i}{2a_2})^4h(\frac{2a_1-\mu_i}{2a_2}).
 \endaligned$$
By Proposition \ref{uniqprop}, one gets under the assumption {\bf A3} that 
$$U_m+\frac{(a_2-r)a_1mr}{2a_2^2}\to N\big(0,\; \frac{\beta \sigma^2}{4}\big)$$
weakly as $a_2\to\infty.$ 
By Lemmas \ref{key-OO}, \ref{key-O-3} and the assumption {\bf A3}, we know 
$$\frac{(a_2-r)^2}{a_2^6}{\rm Var}\big(\sum_{i=1}^m(\mu_i-2a_1)^3\big)=O(\frac{m^7}{m^8})\to 0$$
and $$\frac{(a_2-r)}{24a_2^3}\mathbb{E}\sum_{i=1}^m(\mu_i-2a_1)^3=\frac{(a_2-r)}{24a_2^3}\big(2\beta^2 a_1m^3+o(a_1m^3)\big)\to \frac{\beta^2}{12}xy^3$$
as $a_2\to\infty.$
Consequently, it follows 
$$\frac{(a_2-r)}{24a_2^3}\sum_{i=1}^m(\mu_i-2a_1)^3{\to} \frac{\beta^2}{12}xy^3  $$
in probability as $a_2\to\infty.$
Therefore to prove \eqref{lim-primeO}, it remains to prove 
 \be\lbl{lastprimeO}
\bar{\delta}_m:= (a_2-r)\sum_{i=1}^m(\frac{2a_1-\mu_i}{2a_2})^4h(\frac{2a_1-\mu_i}{2a_2}) \stackrel{p}{\to} 0
\ee
as $a_2\to\infty.$
By Lemma \ref{maxabso}, 
$$\max_{1\le i\le m}|\frac{2a_1-\mu_i}{2a_2}|\to 0$$ in probability $a_2\to\infty.$   
Since $h$ is continuous, $\tau:=\sup_{|x|\leq 1/2}|h(x)|< \infty.$ 
Hence, it follows
\be\lbl{probabO}\aligned 
P(|\bar{\delta}_m|>\epsilon)&= P\Big(|\bar{\delta}_m|>\epsilon,\ \max_{1\leq i \leq m}|\frac{2a_1-\mu_i}{2a_2}|\leq \frac{1}{2}\Big) + P\Big(|\bar{\delta}_m|>\epsilon, \max_{1\leq i \leq m}|\frac{2a_1-\mu_i}{2a_2}|>\frac{1}{2}\Big)\\
& \leq  P\Big(|\bar{\delta}_m|>\epsilon,\ \max_{1\leq i \leq m}|\frac{2a_1-\mu_i}{2a_2}|\leq \frac{1}{2}\Big) + P\Big(\max_{1\leq i \leq m}|\frac{2a_1-\mu_i}{2a_2}|>\frac{1}{2}\Big) \\
\endaligned \ee
as $a_2$ is sufficiently large.
Then under $\max_{1\leq i \leq m}|\frac{2a_1-\mu_i}{2a_2}|\leq \frac{1}{2}$,
$$\aligned 
|\bar{\delta}_m| &\leq  \tau\max_{1\leq i \leq m}|\frac{2a_1-\mu_i}{2a_2}|^2\cdot\frac{(a_2-r)}{4a_2^2}\sum_{i=1}^m(\mu_i-2a_1)^2 \\
&=2\tau \max_{1\leq i \leq m}|\frac{2a_1-\mu_i}{2a_2}|^2 \cdot\frac{(a_2-r)}{8a_2^2}(\sum_{i=1}^m(\mu_i-2a_1)^2-4a_1mr)\\
&\quad +\frac{\tau a_1m^3 r(a_2-r)}{4a_2^4}\max_{1\leq i \leq m}|\frac{2a_1-\mu_i}{m}|^2,
\endaligned 
$$
which converges to zero in probability because $\max_{1\leq i \leq m}|\frac{2a_1-\mu_i}{2a_2}|\to 0$ in probability and $\max_{1\leq i \leq m}|\frac{2a_1-\mu_i}{m}|$ is bounded with probability one
by Lemma \ref{maxabso}, $$\frac{(a_2-r)}{8a_2^2}\big(\sum_{i=1}^m(\mu_i-2a_1)^2-4a_1mr\big)\to N(0, \frac{\beta\sigma^2+\beta^2 xy^3}{4})$$ weakly by Corollary \ref{A3partial} and 
$\frac{a_1m^3 r(a_2-r)}{a_2^4}\to 0$ as $a_2\to\infty$.  This, with \eqref{probabO} and Lemma \ref{maxabso} again concludes \eqref{lastprimeO}.
The proof is close. 
\nprf  

Now we are at the position to post the proof of (ii) of Theorem \ref{bjconv}.  By the relationship \eqref{TV_KL}, it suffices to prove that 
\be\lbl{tv-final}\liminf_{a_2\to\infty}\|\ml{L}(2a\lambda)-\ml{L}(\mu)\|_{\rm TV}>0. \ee  By Lemma 2.14 in \cite{JM2017}, we just need to prove \eqref{tv-final} under assumptions {\bf A2}, {\bf A3} or {\bf A1}. 

\subsubsection{\bf Proof of  (ii) of Theorem \ref{bjconv} under the assumption {\bf A2} or {\bf A3}} 

Review \eqref{mexpre}:
\be\lbl{333}\|\ml{L}(2a\lambda)-\ml{L}(\mu)\|_{\rm TV}=\mathbb{E}|K_m'L_m'(\mu)-1|.\ee Since under {\bf A2} or {\bf A3}, 
$a_1m^3/a^2=a_1m^3/a_2^2+o(1),$
we use Lemma \ref{kmO} to see under {\bf A2} or {\bf A3}, 
$$
\log K_m'=-a_1m+\frac{mr}{2}\log \big(1+\frac {a_1}{a_2}\big)-\frac{\beta^2a_1m^3}{24a_2^2}+o(1)
$$
for $a_2$ large enough.  
This implies \be\lbl{primeof}\aligned  
\log (K_m^{\prime}L_m^{\prime}(\mu))&=\log L_m^{\prime}(\mu)-a_1m+\frac{(a_2-r)a_1mr}{2a_2^2}\\
&-\frac{(a_2-r)a_1m r}{2a_2^2}+\frac{mr}{2}\log \big(1+\frac {a_1}{a_2}\big)-\frac{\beta^2a_1m^3}{24a_2^2}+o(1)\\
\endaligned \ee
for $a_2$ sufficiently large.  
Taylor's  formula allows us to write 
$$\log(1+\frac{a_1}{a_2})=\frac{a_1}{a_2}-\frac{a_1^2}{2a_2^2}+o(\frac{a_1^2}{a_2^2}),$$
which ensures 
$$ \aligned s_m:&=-\frac{(a_2-r)a_1m r}{2a_2^2}+\frac{mr}{2}\log \big(1+\frac {a_1}{a_2}\big)-\frac{\beta^2a_1m^3}{24a_2^2}\\
&=\frac{a_1mr^2}{2a_2^2}-\frac{a_1^2mr}{4a_2^2}-\frac{\beta^2a_1m^3}{24a_2^2}+o(1)\\ 
&=\frac{a_1m (\beta^2 m^2+O(m))}{8a_2^2}-\frac{a_1^2 m(\beta m+2-\beta ) }{8a_2^2}-\frac{\beta^2a_1m^3}{24a_2^2}+o(1)\\
&=\frac{\beta^2 a_1m^3}{12a_2^2}-\frac{\beta a_1^2m^2}{8a_2^2}+\frac{O(a_1m^2+a_1^2m)}{a_2^2}+o(1). \endaligned $$
When the assumption {\bf A2} or {\bf A3} is satisfied,  we know $$\frac{\beta a_1^2m^2}{8a_2^2}\to\frac{\beta\sigma^2}{8} \quad {\rm and} \quad \frac{O(a_1m^2+a_1^2m)}{a_2^2}=O(\frac{1}{m}+\frac{1}{a_1})\to 0$$
as $a_2\to\infty.$   
Then $s_m=\frac{\beta^2 a_1m^3}{12a_2^2}-\frac{\beta\sigma^2}{8}+o(1).$ 
Putting this back to \eqref{primeof}, we have 
$$\log (K_m^{\prime}L_m^{\prime}(\mu))=\log L_m^{\prime}(\mu)-a_1m+\frac{(a_2-r)a_1mr}{2a_2^2}+\frac{\beta^2 a_1m^3}{12a_2^2}-\frac{\beta\sigma^2}{8}+o(1).$$
Since $\frac{\beta^2 a_1m^3}{12a_2^2}\to \frac{\beta^2 xy^3}{12}$ under {\bf A3} or $\frac{a_1m^3}{12a_2^2}\to 0$ under {\bf A2} as $a_2\to\infty,$  
it follows from Lemmas \ref{LemmaLprime} and \ref{centr-O} that 
$$ \log(L_m^{\prime}K_m^{\prime}(\mu))\to N\big(-\frac{\beta \sigma^2}{8}, \frac{\beta \sigma^2}{4}\big)$$
weakly as $a_2\to\infty.$
This implies that $K_m'L_m'(\mu)$ converges weakly to $e^{\xi}$, where $\xi\sim N\big(-\frac{\beta \sigma^2}{8}, \frac{\beta \sigma^2}{4}\big)$. By (\ref{333}) and  the Fatou Lemma, we have 
$$
\liminf_{a_2\to\infty}\|\ml{L}(2a\lambda)-\ml{L}(\mu)\|_{\rm TV}\geq \mathbb{E}|e^{\xi}-1|>0.
$$
The proof is finished. 
$\hfill\square $
\subsubsection{\bf Proof of (ii) of Theorem \ref{bjconv} under the assumption {\bf A1}}  For this particular case, we know $r=m=1$ and $\frac{a_1}{a_2}\to \sigma\in(0,1).$ Therefore Lemma \ref{K} tells 
$$\log K_1^{\prime}=-a_1+\frac12\log(1+\sigma)+o(1)$$ 
and with the help of  Taylor's expansion $$\log L_1^{\prime}(\mu)=a_1+\frac{\mu-2a_1}{2a_2}-\frac{a_2-1}{8a_2^2}(\mu-2a_1)^2+(a_2-1)\frac{(2a_1-\mu)^3}{a_2^3}h(\frac{\mu-2a_1}{2a_2})$$
with $h$ a continuous function on $(-1, \infty).$  Here $\mu$ has density function $f_{\beta, a_1}$ with $m=r=1.$
Then 
\be\lbl{exp-c}\log (K_1^{\prime} L_1^{\prime}(\mu))=\frac12\log(1+\sigma)+\frac{\mu-2a_1}{2a_2}-\frac{a_2-1}{8a_2^2}(\mu-2a_1)^2+(a_2-1)\frac{(2a_1-\mu)^3}{a_2^3}h(\frac{\mu-2a_1}{2a_2})+o(1).\ee
Examining the form $f_{\beta, a_1}$ in this particular case, 
we see $\mu\sim\chi^2_{2a_1}.$  That means we could rewrite $\mu-2a_1$ as $\mu-2a_1=\sum_{i=1}^{2a_1}(\xi_i^2-1)$ with $\xi_i\sim N(0, 1)$ for $1\le i\le 2a_1$ and $(\xi_i)_{1\le i\le 2a_1}$ are mutually independent.
Since $\mathbb{E}\xi_i^2=1$ and ${\rm Var}(\xi_i^2-1)=2,$ by Lindeberg-L\'evy central limit Theorem, 
 we see 
$$\frac{\mu-2a_1}{2\sqrt{a_1}}\to N(0, 1)$$ weakly as $a_2\to\infty.$ 
Consequently, $$\frac{\mu-2a_1}{2a_2}=\frac{\mu-2a_1}{2\sqrt{a_1}}\frac{\sqrt{a_1}}{a_2}\to 0 \quad  \text{and} \quad \frac{(\mu-2a_1)^3}{a_1^2}=(\frac{\mu-2a_1}{\sqrt{a_1}})^3\frac{1}{\sqrt{a_1}}$$ in probability as $a_2\to\infty$ and then $\frac{(2a_1-3)^3}{a_2^2}h(\frac{\mu-2a_1}{2a_2})$ tends to $0$ in probability as $a_2\to\infty.$ 
For the term $\frac{a_2-1}{8a_2^2}(\mu-2a_1)^2,$ similarly we have
$$\frac{a_2-1}{8a_2^2}(\mu-2a_1)^2=(\frac{\mu-2a_1}{2\sqrt{a_1}})^2\frac{a_1(a_2-1)}{2a_2^2}\to \frac{\sigma}{2} \chi^2_1$$
weakly as $a_2\to\infty.$
Putting all these limits into \eqref{exp-c}, we know 
$$\frac{f_{\beta, a_1, a_2}}{f_{\beta, a_1}}=e^{\log(K_1^{\prime}L_1^{\prime}(\mu))}\to\sqrt{1+\sigma}\exp\{-\frac{\sigma}{2}\chi^2_1\} $$
weakly as $a_2\to\infty.$
By \eqref{333} and the Fatou Lemma, 
$$\aligned \liminf_{a_2\to\infty}\|{\ml L}(2a\lambda)-{\ml L}(\mu)\|_{\rm TV}&\ge \mathbb{E}|\sqrt{1+\sigma}e^{-\frac{\sigma}{2}\chi^2_1}-1|
>0.
\endaligned $$
Finally the whole proof is  close now. \hfill$\square$
  
\section{Appendix}

\bprf[\bf Proof of Lemma \ref{key-OO}.]  Review that $x_i^2\sim \chi^2_{(2a_1-\beta (i-1))}$ and $y_i^2\sim\chi^2_{\beta(m-(i-1))}.$

Setting $b_i:=\beta m-2\beta(i-1)$  for $2\le i\le m$ and $b_1=0,$  one gets $$z_i:=x_i^2+y_i^2\sim \chi^2_{2a_1+b_i}$$ 
for $1\le i\le m$ with the convention $y_1=0.$

It is easy to see \be\lbl{zero}\aligned &\sum_{i=1}^{m}b_i=\beta\sum_{i=1}^{m-1}(m-2i)=0;\\
&\sum_{i=1}^m b_i^2=\beta^2\sum_{i=1}^{m-1}(m-2i)^2=\frac{\beta^2 m(m-1)(m-2)}{3}; \\
&\sum_{i=1}^m b_i^3=\beta^3\sum_{i=1}^{m-1}(m-2i)^3=0.\\
\endaligned \ee
Based on the Dumitriu and Edelman characterization, one gets 
\be\lbl{form-single}\sum_{i=1}^m \mu_i={\rm tr}({\bf AA}^{\prime})=\sum_{i=1}^{m}(x_{i}^2+y^2_{i})=\sum_{i=1}^m z_i. 
\ee
Therefore by Lemma \ref{chi} and \eqref{zero},  we have 
\be\lbl{exp}\aligned \mathbb{E}\sum_{i=1}^{m} \mu_i=\sum_{i=1}^{m}\mathbb{E}z_i
=\sum_{i=1}^{m}(2a_1+b_i)=2a_1m
\endaligned \ee
and consequently 
\be\lbl{cov} \aligned &{\rm Var}\bigg(\sum_{i=1}^m\mu_i\bigg)=\sum_{i=1}^{m}{\rm Var}(z_i)=2\sum_{i=1}^m\mathbb{E}z_i=4a_1m.
\endaligned \ee
By Dumitriu and Edelman's characterization again, $\mu_1-2a_1, \mu_2-2a_1, \cdots, \mu_m-2a_1$ are the eigenvalues of the matrix ${\bf AA}'-2a_1 {\bf I}_{m}.$  
Therefore  
\be\lbl{sumofsquare}\aligned \sum_{i=1}^m(\mu_i-2a_1)^2={\rm tr}\big(({\bf AA}^{\prime}-2a_1 {\bf I}_m)^2\big)
=\sum_{i=1}^{m}(z_i-2a_1)^2+2\sum_{i=2}^{m}x_{i-1}^2y_i^2.\\
\endaligned 
\ee 
Remembering $\mathbb{E}z_i=b_i+2a_1$ for $1\le i\le m,$ we have  from Lemma \ref{chi}
$$\mathbb{E}(z_{i}-2a_1)^2=\mathbb{E}(z_{i}-\mathbb{E}z_i+b_i)^2={\rm Var}(z_i)+b_i^2=4a_1+2b_i+b_i^2
$$
for $1\le i\le m.$ This and 
\eqref{zero} bring us 
$$\sum_{i=1}^m\mathbb{E}(z_{i}-2a_1)^2=4a_1m+\frac{\beta^2 m(m-1)(m-2)}{3}.$$
By independence, it is clear that  
$$\aligned 2\sum_{i=2}^m\mathbb{E}(x_{i-1}^2 y_i^2)&=2\beta \sum_{i=2}^{m}(2a_1-\beta(i-2))(m-(i-1))\\
&=2\beta a_1m(m-1)-\frac{\beta^2m(m-1)(m-2)}{3}.\endaligned $$
Then we have   
$$\mathbb{E}\sum_{i=1}^m (\mu_{i}-2a_1)^2=2a_1m(2+\beta(m-1))=4a_1mr. $$ 
Now we investigate the third expression. It follows from \eqref{sumofsquare} that 
\be\lbl{sum22} \aligned {\rm Var}(\sum_{i=1}^m(\mu_i-2a_1)^2)&=\sum_{i=1}^m {\rm Var}((z_i-2a_1)^2)+4\sum_{i=2}^m {\rm Var}(x_{i-1}^2 y_i^2)\\
&\quad +4{\rm Cov}\big(\sum_{i=1}^m (z_i-2a_1)^2, \sum_{i=2}^mx_{i-1}^2 y_i^2\big). \endaligned \ee
Next we examine one by one the three terms in \eqref{sum22}.  
 Using $\mathbb{E}z_i=b_i+2a_1$  and Lemma \ref{chi} again, we know 
 \be\lbl{varz}\aligned {\rm Var}((z_i-2a_1)^2)&={\rm Var}\big((z_i-\mathbb{E}z_i+b_i)^2\big)\\
&={\rm Var}\big((z_i-\mathbb{E}z_i)^2)+4b_i^2{\rm Var}(z_i)+ 4b_i{\rm Cov}((z_i-\mathbb{E}z_i)^2, z_i-\mathbb{E}z_i)\big)\\
&=8\mathbb{E}z_i(\mathbb{E}z_i+6)+8b_i^2\mathbb{E}z_i+32b_i\mathbb{E}z_i\\
&=8\big(4a_1^2+12a_1+(12a_1+6)b_i+(2a_1+5)b_i^2+b_i^3\big).
\endaligned \ee
By this expression and \eqref{zero}, we have  
\be\lbl{var-expression} \sum_{i=1}^{m} {\rm Var}((z_i-2a_1)^2)=32a_1(a_1+3)m+\frac{8\beta^2}{3}(2a_1+5)m(m-1)(m-2).
 \ee
For the second term $2\sum_{i=1}^{m-1}x_{i}^2y_{i+1}^2$ in \eqref{sum22}, by independence, we have 
\be\lbl{varinter}\aligned {\rm Var} (x_{i}^2y_{i+1}^2)&=\mathbb{E}[x_{i}^4y_{i+1}^4]-(\mathbb{E}[x_{i}^2y_{i+1}^2])^2\\
 &=\big((\mathbb{E}x_i^2)^2+2\mathbb{E}x_i^2\big)\big((\mathbb{E}y_{i+1}^2)^2+2\mathbb{E}y_{i+1}^2\big)-(\mathbb{E}x_i^2)^2(\mathbb{E}y_{i+1}^2)^2\\
 &=2\mathbb{E}x_i^2\mathbb{E}y_{i+1}^2(\mathbb{E}x_i^2+\mathbb{E}y_{i+1}^2+2)\\
 &=2(2a_1+\beta-\beta i)\beta (m-i)(2a_1+\beta+2+\beta (m-2i)).
\endaligned \ee 
Therefore it follows from careful calculation that  
 \be\lbl{inter-expression} \aligned  {\rm Var} \big(\sum_{i=1}^{m-1}x_{i}^2y_{i+1}^2\big)
 &=4\beta a_1^2m(m-1)+4\beta a_1m(m-1)+2\beta^2m(m-1)(a_1-\frac{m-2}3) . \endaligned \ee
Now we work on the last term in \eqref{sum22}.  On the one hand,  we have 
\be\lbl{up-down1}\aligned  {\rm Cov}((z_i-2a_1)^2, y_i^2)&={\rm Cov}\big((y_i^2-\mathbb{E}y_i^2+x_i^2-\mathbb{E}x_i^2+b_i)^2, y_i^2-\mathbb{E} y_i^2\big)\\
&=\mathbb{E}(y_i^2-\mathbb{E}y_i^2)^3+2b_i{\rm Var}(y_i^2)\\
&=4(2+b_i) \mathbb{E}y_i^2\endaligned \ee
for $1\le i\le m.$ Here the last equality is guaranteed again by Lemma \ref{chi} and the second one is true since $x_i^2-\mathbb{E}x_i^2+b_i$ is independent of $y_i.$ 
On the other hand, similarly we have 
\be\lbl{up-down2} {\rm Cov}((z_i-2a_1)^2, x_i^2)=4(2+b_i) \mathbb{E}x_i^2 \ee
for $1\le i\le m.$
Therefore by independence and \eqref{up-down1} and \eqref{up-down2},  we have 
$$\aligned & \quad\quad {\rm Cov}\bigg(\sum_{i=1}^m (z_i-2a_1)^2, \sum_{i=2}^mx_{i-1}^2 y_i^2\bigg)\\
&=\sum_{i=1}^{m-1}\mathbb{E} y_{i+1}^2{\rm Cov}\big((z_i-2a_1)^2, x_{i}^2\big)+\sum_{i=2}^m \mathbb{E}x_{i-1}^2{\rm Cov}\big((z_{i}-2a_1)^2, y_{i}^2\big)\\\
&=4\sum_{i=1}^{m-1} \mathbb{E} y_{i+1}^2 \mathbb{E}x_{i}^2(4+b_i+b_{i+1}).\\
\endaligned $$ 
Thereby with simple algebra, we have   
\be\lbl{cov-expression}\aligned {\rm Cov}\bigg(\sum_{i=1}^m (z_i-2a_1)^2, \sum_{i=2}^mx_{i-1}^2 y_i^2\bigg)&=\frac{8\beta^2}{3}(a_1-1)m(m-1)(m-2)
+16\beta a_1 m(m-1).\endaligned \ee
Plugging \eqref{var-expression},  \eqref{inter-expression} and \eqref{cov-expression} into \eqref{sum22}, we finally have 
$$\aligned {\rm Var}(\sum_{i=1}^m (\mu_i-2a_1)^2)&=16\beta a_1^2m(m-1)+16 \beta^2 a_1m(m-1)(m-2)+8\beta^2a_1 m(m-1)\\
&\quad+80\beta a_1m(m-1)+32a_1m(a_1+3). 
\endaligned $$ 
Now we prove the expression for covariance.  Similarly since all the random variables involved are independent, 
we have  by \eqref{form-single} and \eqref{sumofsquare},
$$\aligned
{\rm Cov}(\sum_{i=1}^m (\mu_i-2a_1)^2, \sum_{i=1}^m (\mu_i-2a_1))
&=\sum_{i=1}^m{\rm Cov}\big((z_i-2a_1)^2, x_i^2+y_i^2\big)\\
&+2\sum_{i=1}^{m-1}{\rm Cov}(x_{i}^2y_{i+1}^2, x_{i}^2+y_{i+1}^2). \\
\endaligned $$
Then \eqref{up-down1}, \eqref{up-down2},  the independence of $x_i$ and $y_j$ and Lemma \ref{chi} show that 
$$
\aligned &{\rm Cov}(\sum_{i=1}^m (\mu_i-2a_1)^2, \sum_{i=1}^m (\mu_i-2a_1))\\
&=\sum_{i=1}^m(8+4b_i)\mathbb{E}z_i+8\sum_{i=1}^{m-1}\mathbb{E}y_{i+1}^2\mathbb{E}x_{i}^2\\
&=4\sum_{i=1}^m(b_i+2a_1)(b_i+2)+8\sum_{i=1}^{m-1}\beta(m-i)(2a_1-\beta (i-1)).
\endaligned$$
With simple calculus on the sum, we get from \eqref{zero}
$${\rm Cov}(\sum_{i=1}^m (\mu_i-2a_1)^2, \sum_{i=1}^m (\mu_i-2a_1))=8\beta a_1m^2+8a_1m(2-\beta)=16a_1mr.$$
It remains to prove the last expression.  By the property of the random matrix ${\bf A},$ it is not hard to verify that 
\be\lbl{formula3}\aligned \sum_{i=1}^m(\mu_i-2a_1)^3&={\rm tr}\big( ({\bf AA'}-2a_1{\bf I}_m)^3 \big)\\
&=\sum_{i=1}^m(z_i-2a_1)^3+3\sum_{i=1}^{m-1}x_i^2y_{i+1}^2(z_i+z_{i+1}-4a_1). 
\endaligned 
\ee
For the first term. It is ready to check that 
$$\aligned \mathbb{E}(z_i-\mathbb{E}z_i+b_i)^3&=\mathbb{E}(z_i-\mathbb{E}z_i)^3+3b_i{\rm Var}(z_i)+b_i^3\\
&=8\mathbb{E}z_i+6b_i\mathbb{E}z_i+b_i^3\\
&=8(2a_1+b_i)+6b_i(2a_1+b_i)+b_i^3
\endaligned $$ 
for $1\le i\le m.$ Thereby with the help of \eqref{zero}, we have   
$$\mathbb{E}\sum_{i=1}^m(z_i-2a_1)^3=16a_1m+2\beta^2m(m-1)(m-2).$$
Also by Lemma \ref{chi}, it follows  
$$\aligned \mathbb{E}[x_i^2y_{i+1}^2(z_i+z_{i+1}-4a_1)]&=\mathbb{E}[x_i^2y_{i+1}^2(x_i^2+y_i^2+x_{i+1}^2+y_{i+1}^2-4a_1)]\\
&=\mathbb{E}x_i^2\mathbb{E}y_{i+1}^2 \big(\mathbb{E}(x_i^2+y_i^2+2)+\mathbb{E}(x_{i+1}^2+y_{i+1}^2+2)-4a_1\big)\\
&=\mathbb{E}x_i^2\mathbb{E}y_{i+1}^2 (b_i+b_{i+1}+4)
\endaligned $$
for any $1\le i\le m-1.$ 
 Therefore by \eqref{cov-expression}, we have 
 $$\aligned 3\sum_{i=1}^{m-1}\mathbb{E}x_i^2y_{i+1}^2(z_i+z_{i+1}-4a_1)&=
2\beta^2(a_1-1)m(m-1)(m-2)
+12\beta a_1 m(m-1). \endaligned $$ 
Consequently, we have  
$$\mathbb{E}\sum_{i=1}^m (\mu_i-2a_1)^3=2\beta^2a_1m(m-1)(m-2)+12\beta a_1m(m-1)+16a_1m.
$$
The proof is complete now. 
\nprf 

\bprf[\bf Proof of Lemma \ref{key-O-3}.] 
By \eqref{formula3} and the property of variance, we have 
\be\lbl{zeroorder}
\aligned 
{\rm Var}(\sum_{i=1}^m (\mu_i-2a_1)^3)\le2\sum_{i=1}^m {\rm Var}((z_i-2a_1)^3)+18(m-1)\sum_{i=1}^{m-1}{\rm Var}(x_i^2y_{i+1}^2(z_i+z_{i+1}-4a_1)).
\endaligned \ee
On the one hand, we have 
\be\lbl{firstorder}\aligned {\rm Var}((z_i-2a_1)^3)&\le \mathbb{E}(z_i-\mathbb{E}z_i+b_i)^6\\
&\le 32 \mathbb{E}(z_i-\mathbb{E}z_i)^6+32 b_i^6\\
&=O(m^6),
\endaligned \ee
where the last equality holds since $a_1=O(m)$ and $b_i=O(m)$ for all $1\le i\le m.$  
On the other hand, it holds 
$$\aligned {\rm Var}(x_i^2y_{i+1}^2(z_i+z_{i+1}-4a_1))&\le 3{\rm Var}\big(x_i^2y_{i+1}^2(z_i-\mathbb{E}z_i)\big)+3{\rm Var}\big(x_i^2y_{i+1}^2(z_{i+1}-\mathbb{E}z_{i+1})\big)\\
&\quad +3(b_i+b_{i+1})^2{\rm Var}(x_i^2y_{i+1}^2)\\
\endaligned $$
for $1\le i\le m-1.$ 
By \eqref{varinter}, $$3(b_i+b_{i+1})^2{\rm Var}(x_i^2y_{i+1}^2)=O(m^5).$$ 
Obviously 
$$\aligned {\rm Var}\big(x_i^2y_{i+1}^2(z_i-\mathbb{E}z_i)\big)&\le\mathbb{E}x_i^4y_{i+1}^4(z_i-\mathbb{E}z_i)^2 \\
& =\mathbb{E}y_{i+1}^4\mathbb{E}x_i^4(z_i-\mathbb{E}z_i)^2.
\endaligned $$
By the independence of $x_i$ and $y_i,$ we have 
$$\aligned \mathbb{E}x_i^4(z_i-\mathbb{E}z_i)^2&=\mathbb{E}x_i^4(y_i^2-\mathbb{E}y_i^2)^2+\mathbb{E}x_i^4(x_i^2-\mathbb{E}x_i^2)^2\\
&=\mathbb{E}x_i^4 {\rm Var}(y_i^2)+\mathbb{E}(x_i^2-\mathbb{E}x_i^2)^4+(\mathbb{E}x_i^2)^2{\rm Var}(x_i^2)+2 \mathbb{E}x_i^2 \mathbb{E}(x_i^2-\mathbb{E}x_i^2)^3\\
&=2\mathbb{E}x_i^4 \mathbb{E}y_i^2+12\mathbb{E}x_i^2(\mathbb{E}x_i^2+4)+2(\mathbb{E}x_i^2)^3+16(\mathbb{E}x_i^2)^2\\
&=O(m^3),
\endaligned $$
where the third equality is guaranteed by Lemma \ref{chi}.  This ensures 
$${\rm Var}\big(x_i^2y_{i+1}^2(z_i-\mathbb{E}z_i)\big)=O(m^5).$$
Similarly we have 
$${\rm Var}\big(x_i^2y_{i+1}^2(z_{i+1}-\mathbb{E}z_{i+1})\big)=O(m^5).$$ 
Therefore 
\be\lbl{secondorder}{\rm Var}(x_i^2y_{i+1}^2(z_i+z_{i+1}-4a_1))=O(m^5).\ee 
Hence combining \eqref{zeroorder}, \eqref{firstorder} and \eqref{secondorder}, we know 
$${\rm Var}(\sum_{i=1}^m (\mu_i-2a_1)^3)=O(m^7).$$ 
The proof is complete.  
\nprf

\end{document}